
\input amstex
\documentstyle{amsppt}
\NoRunningHeads
 \magnification\magstep1
 \TagsOnLeft
\define\n{\noindent}

\overfullrule = 0pt

\define\vp{\varepsilon}
\define\ent{{\text{\rm Z}\mkern-5.5mu\text{\rm Z}}}
 \overfullrule = 0pt

\define\reel{\text{\rm I}\!\text{\rm R}}

\define\comp{\;{}^{ {}_\vert }\!\!\!\text{\rm C}}
\define\N#1{\left\Vert#1\right\Vert}
\define\nat{ { \text{\rm I}\!\text{\rm N}} }
\define\ee{ { \text{\rm I}\!\text{\rm E}}}
\define\T{{\bold T}}

\define\ie{{\it i.e.\/}\ }
\define\cf{{\it cf.\/}\ }

\topmatter
\title The similarity degree of an operator algebra\endtitle
\author by Gilles Pisier\endauthor
\address Texas A\&M University,
College Station, TX 77843, U. S. A.\
and
Universit\'e Paris VI,
Equipe d'Analyse, Case 186, 75252,
 Paris Cedex 05, France\endaddress
\thanks Partially supported by the N.S.F.\endthanks
\thanks Revised April 97, submitted to St.\ Petersburg Math J.\endthanks
\toc
\specialhead{} Plan\endspecialhead
\head \S 0. Introduction\endhead
\head \S 1. Enveloping operator algebras. Preliminary results\endhead
\head \S 2. Main results\endhead
\head \S 3. Groups \endhead
\head \S 4. Operator algebras \endhead
\head \S 5. Uniform algebras   \endhead
\head \S 6. $C^*$-algebras \endhead
\head \S 7.  The Blecher-Paulsen factorization \endhead
\head \S 8.  Banach algebras \endhead
\endtoc
\abstract
 Let $A$ be a unital operator algebra.
Let us assume that every {\it bounded\/} unital 
homomorphism $u\colon \ A\to B(H)$ is similar to
a {\it contractive\/} one.
 Let 
$\text{\rm Sim}(u) = \inf\{\|S\|\, \|S^{-1}\|\}$
where the infimum runs over all invertible operators $S\colon \ H\to H$ such 
that the
``conjugate'' homomorphism $a\to S^{-1}u(a)S$ is contractive. Now for all
$c>1$, let
$\Phi(c) = \sup\text{\rm Sim}(u)$
where the supremum runs over all unital homomorphism $u\colon\ A\to B(H)$
with $\|u\|\le c$.
  Then,  there is
$\alpha\ge 0$ such that for some constant $K$ we have:
$$\Phi(c) \le Kc^\alpha.\leqno (*)\qquad \forall c>1$$
Moreover, the smallest $\alpha$ for which this holds is an integer, denoted
by $d(A)$ (called the similarity degree of $A$) and $(*)$ still holds for some 
$K$ when $\alpha=d(A)$.
Among the applications of these  results,   we give new
characterizations of proper uniform algebras on one hand, and of nuclear
$C^*$-algebras on the other. Moreover, we obtain a characterization of amenable
groups  which answers (at least partially) a question on group
representations going back to a 1950 paper of Dixmier.
\endabstract
\endtopmatter

\head\S 0. Introduction\endhead

Consider a unital operator algebra $A$ (i.e.\ a subalgebra of $B(H)$,
containing $I$, not assumed self-adjoint). 
We are interested in the following ``similarity property" of $A$:

For any {\it bounded} unital homomorphism $u\colon \ A\to B(H)$, 
there is an invertible
operator \break 
\indent $S: H\to H$ (= a similarity)
such that $x\to S^{-1} u(x) S$ is {\it contractive}. 

\n In other words,
 every {\it bounded} unital homomorphism on $A$ is similar
to a {\it contractive} one.

\n Let 
$\text{\rm Sim}(u) = \inf\{\|S\|\, \|S^{-1}\|\}$
where the infimum runs over all invertible
 operators $S\colon \ H\to H$ such that the
``conjugate'' homomorphism $a\to S^{-1}u(a)S$ is 
contractive. Now for all
$c>1$, let
$$\Phi(c) = \sup\text{\rm Sim}(u)$$
where the supremum runs over all unital 
homomorphism $u\colon\ A\to B(H)$
with $\|u\|\le c$.
 Assume that the above similarity property holds.
 Then it is easy
to show that $\Phi(c)$ is finite for all $c>1$. 
Our first observation,
simple but crucial, will be that necessarily 
$\Phi(c)$ has {\it polynomial
growth\/}, i.e.\ there is a number $\alpha\ge 0$ and a constant $K$   such
that 
$$\forall c>1 \quad \Phi(c)\le K c^\alpha ,\tag0.1$$
 equivalently:  any bounded unital homomorphism
 $u\colon \ A\to B(H)$ satisfies
$\text{\rm Sim}(u) \le K\|u\|^\alpha.$
Let $d(A)$ be the infimum of the numbers $\alpha\ge 0$ for which
(0.1) holds for some constant  $K$.

Our second observation (which lies a bit deeper) is that $d(A)$ is an
integer, i.e.\ we have
$$d(A) \in \{0,1,2,3,\ldots\},$$
and moreover there is a constant $K$ such that (0.1) holds for $\alpha =
d(A)$.

We call $d(A)$ the similarity degree of the operator
 algebra $A$. If the similarity property fails, then
we set $d(A) =\infty$.

By a result due to Paulsen ([Pa4]), the similarity property is closely related 
to
the notion of  complete  boundedness,
 for  which we refer to  [Pa1]. 
To decribe this connection, we will consider
the following property (C) of an operator algebra $A$:

(C) Every contractive unital homomorphism
$u\colon \ A\to B(H)$ is completely bounded. 

  It is easy to see that this holds  for all
$C^*$-algebras and   for several examples
 of uniform algebras
 (such as the disc and the bidisc algebras).

Under this assumption, (see [Pa4]) a unital homomorphism
 $u\colon \ A\to B(H)$ is similar to a contractive one iff
it is completely bounded ($c.b.$ in short). 

Let $\Cal K$ be the $C^*$-algebra of compact
 operators on $\ell_2$, let
$C_0 \subset \Cal K$ be the subspace of diagonal
operators and let
$\Cal K \otimes_{\text{\rm min}} A$ be the minimal (= spatial) tensor product.
Under the above assumption (C) on $A$, we will show (see Theorem 4.2) that 
$d(A)$ is
the smallest integer $d$ with the following property:\ there is a constant
$K$ such that any $x$ in the unit ball of $\Cal K\otimes_{\text{\rm min}} A$
can be written as a product of the form
$$x = \alpha_0 D_1\alpha_1D_2\ldots
D_d\alpha_{d}$$ with $\alpha_i \in \Cal
K\otimes 1$ and $D_i \in C_0 \otimes_{\text{\rm min}}
A$ such that
$$\prod^{d}_{i=0} \|\alpha_i\| \prod^d_{i=1} \|D_i\| \le K.$$
Thus, $d(A)$ appears as the minimal ``length" necessary to
express  any element of the unit ball of 
$\Cal K \otimes_{\text{\rm min}} A$
as a alternated product as above with $2d+1$ factors
(with a good control of the norms of the factors).

More generally, if $\Cal A$ is merely a Banach algebra
with unit, we may consider it as embedded as a dense
unital subalgebra into its enveloping 
unital operator algebra $\tilde A$. The morphism
 $\Cal A\subset \tilde A$ is characterized by the property that 
a unital
homomorphism $v\colon \ \Cal A\to B(H)$ is contractive
 (\ie has norm equal to 1) iff it extends  
 to a completely contractive 
homomorphism $\tilde v\colon \tilde A\to B(H)$.
In particular, $\tilde A$ satisfies (C).

\n In this situation, let us assume that every bounded unital homomorphism
$u\colon \ \Cal A\to B(H)$ extends to a completely bounded unital homomorphism
$\tilde u\colon \ \tilde A\to B(H)$. We define
$$d(\Cal A) = \inf\{\alpha\ge 1\}\tag0.2$$
where the infimum runs over all $\alpha\ge 1$ such that for some $K$ we
have for all bounded unital homomorphisms $u\colon \ \Cal A\to B(H)$
$$\|\tilde u\|_{cb} \le K\|u\|^\alpha.\tag0.3$$
If there is no such $\alpha$, then we set by convention $d(\Cal A) = \infty$.
Then again the same observations are valid:\medskip
\roster
\item"{(i)}" $d(\Cal A)<\infty$
\item"{(ii)}" $d(\Cal A)$ is an integer and the infimum is attained in
(0.2).
\endroster\medskip

\n An interesting example of this situation is given by group algebras (or
semi-group algebras).

Let $G$ be a discrete group (resp.\ semi-group with
unit). Let $\Cal A$ be the group (resp.\ semi-group) algebra of
$G$ i.e.\ $\Cal A=\ell_1(G)$ equipped with convolution.
In the group case, $\tilde A$ coincides with the (full)
$C^*$-algebra of $G$, denoted by $C^*(G)$.  Let
$g\to \delta_g$ be the natural mapping from $G$ into
$\ell_1(G)$ (i.e.\ $\delta_g(s)  = 1$ iff $s=g$). Let
$u\colon \ \ell_1(G)\to B(H)$ be a linear map and let
$\pi(g) = u(\delta_g)$. Clearly $u$ is a bounded unital
homomorphism iff $\pi$ is a uniformly bounded group
(resp.\ semi-group) representation. Moreover if we define
$$|\pi| = \sup_{g\in G} \|\pi(g)\|$$ we have obviously
$$|\pi| = \|u\|.$$ Conversely, any bounded
representation $\pi\colon \ G\to B(H)$ extends uniquely to
a bounded linear homomorphism $u\colon \ \ell_1(G)\to
B(H)$.

In this setting, we will write $d(G)$ instead of $d(\Cal A)$.

 \n We can show (see Theorem 3.2
and Corollary 3.4 below) that $d(G)=1$ iff $G$ is finite
and $d(G)=2$ iff $G$ is amenable and infinite.

 This result
gives some information on the ``similarity problem'' for uniformly bounded
group representations. Namely, we can prove

\proclaim{Theorem 0.1} Let $G$ be a discrete group. The following are
equivalent:\medskip
\roster
\item"{\text{\rm (i)}}" $G$ is amenable.
\item"{\text{\rm (ii)}}" There is a constant $K$ and $\alpha<3$ such that for 
any $H$
and for any uniformly bounded group representation $\pi\colon \ G\to B(H)$
there is an invertible operator $S\colon \ H\to H$ (called ``a
similarity'') with
$$\|S^{-1}\|\, \|S\| \le K|\pi|^\alpha$$
and such that
$g\to S^{-1}\pi(g)S$
is a unitary representation of $G$.
\item"{\rm (iii)}" Same as {\rm (ii)} with $K=1$ and $\alpha=2$.\medskip
\endroster
\endproclaim

\n Note:\ A uniformly
bounded representation $\pi\colon  \ G\to B(H)$ is called unitarizable if
there is an invertible $S\colon \ H\to H$ such that $S^{-1}\pi(\cdot)S$ is
a unitary representation. The implication (i)~$\Rightarrow$~(iii) is a
classical fact proved in 1950 by Dixmier [Di], following earlier work
by Sz.- Nagy [SN] for $G=\ent$.
 At that time, 
there were no known example of uniformly bounded  
non-unitarizable representation.
  The first example of this phenomenon
was given in 1955 by Ehrenpreis and Mautner [EM] 
on the group $SL_2(\reel)$
 (\cf also     [KS]). See Cowling's notes [Co]
for more information on the Lie group case.
 Later on, many constructions were
given on non-commutative free groups (or on any discrete
group containing
 a non-commutative free group as a subgroup). See 
for example the references in [MPSZ] and [BF2].
See also [P1, Chapter2].

\n In the same paper,
 Dixmier asks whether amenable groups are
the only groups $G$ on which every  uniformly bounded  
representation $\pi$ is unitarizable. This  remains an open question.
 Our result shows that if
one incorporates in Dixmier's question the fact that the similarity 
$S$ can
be found with $\|S\|\, \|S^{-1}\| \le |\pi|^2$, then the
answer is affirmative.

It  seems conceivable that $d(G)<\infty \Rightarrow d(G)\le 2$ automatically, 
but at
the time of this writing we
 have not been able to prove this, and we are now more
 inclined to believe
(in analogy with Corollary 6.2)
 that there are examples of discrete groups $G$
with $2<d(G)<\infty$. Note 
that these would be non-amenable
groups not containing ${\bold F}_2$, the free group  on two generators.
While such examples are known to exist [O1-2],
they still appear difficult to understand (see for
example  the exposition in [Pat]).

Recently, we proved ([P5])
 that when $A$ is the disc algebra
we have $d(A)=\infty$, thus solving the ``Halmos problem" 
on polynomially bounded operators.
Of course, this also holds for the polydisc algebra,
 the ball algebra or for
any uniform algebra  admitting 
a quotient algebra (unitally) isometric to (or
completely isomorphic to)  the disc algebra.
It is conceivable that $d(A)=\infty$ for 
any proper uniform algebra, however
at this point we are only able to show  the following (see Theorem 5.1).

\proclaim{Theorem 0.2}  Let $K$ be a compact set. Let
$A\subset C(K)$ be a uniform algebra (\ie a closed unital
subalgebra which separates the points of $K$).
Then $A=C(K)$ iff $d(A)\le 2$ and $A$  satisfies (C).
\endproclaim

We now turn to $C^*$-algebras. Unfortunately,
at this time we are  unable
   to produce examples
of $C^*$-algebras $A$ for which $d(A)$ takes
arbitrarily large finite values (or one for which
the degree
is infinite). This
 would solve (negatively) a well
known open problem, due to Kadison [Ka] (see [P1]).
 We conjecture that there is a $C^*$-algebra $A$
(probably the
reduced $C^*$-algebra 
of the free group on   infinitely many generators) with
$d(A)=\infty$. 
Unfortunately we only are able to produce examples
of $C^*$-algebras $A$ with $d(A)$ equal to either 1 
(finite dimensional case), 2 (nuclear case), and 3
 ($B(\ell_2)$). 

We give a result (Theorem 6.1) 
which is very close to proving
that, for a $C^*$-algebra,  $d(A)\le 2$
 implies $A$ nuclear. Indeed, it is known (see [CE]) that $A$ is nuclear
iff, for any $*$-representation $\pi\colon\ A\to B(H)$,
the von Neumann algebra generated by $\pi$ is injective.
What we can prove is the following (see Theorem 6.1).

\proclaim{Theorem 0.3} Let $A$ be a unital $C^*$-algebra
such that $d(A)\le 2$. Then,
 whenever a $*$-representation $\pi\colon\ A\to B(H)$
 generates a semi-finite von Neumann algebra,  that
von Neumann algebra is injective.
\endproclaim

We are convinced that the similarity degree $d(A)$ can take arbitrary integer 
values
when $A$ runs over all possible (non self-adjoint) operator algebras,
but again we have not been able to verify this yet. 
However, in the more general
 framework of ``similarity settings" considered below,
 it is easy to exhibit
examples realizing  any possible integral value of the degree,
see Remark 3.6.

The present investigation was considerably influenced by several sources
 [Pel, B1, BRS, BP2] which I would like to acknowledge here:\medskip

\n 1)~~Peller's paper [Pel] contains a discussion (partly based on some
ideas of A.~Davie [Da]) of the space of coefficients of representations of a
$Q$-algebra (with some consequences for operator
algebras). In view of the recent characterizations in
[BRS] and [B1] of operator algebras which use the Haagerup
tensor product, it was natural to try to transpose these
ideas from [Pel] to the ``new'' category of operator
algebras with c.b.\ maps as its morphisms. This is the
content of section~1 below. 

\n 2)~~Blecher and Paulsen's paper [BP2] contains several striking
factorization theorems for elements in the maximal tensor products of
various operator algebras, analogous to the factorization of polynomials
into products of polynomials of degree 1. 
Their factorization is into infinitely (or at least unboundedly) many
matricial factors (see \S 7). It was natural to wonder in which case the number 
of
factors could be bounded by a fixed integer. This is what lead to the
central notion of this paper:\ the ``similarity degree'' (see Theorem 4.2).
\medskip

We refer the reader to the books [Pa1] and [P1] for the precise
definitions of all the undefined terminology
 that we will use, and to [KaR] for operator algebras in general. We recall
only that an ``operator space'' is a closed subspace $E\subset B(H)$ of
the $C^*$-algebra of all bounded operators on a Hilbert space $H$. We will
use freely the notion of a completely bounded (in short c.b.) map $u\colon
\ E_1\to E_2$ between two operator spaces, as defined e.g.\ in [Pa1]. We
denote by $\|u\|_{cb}$ the corresponding norm and by $cb(E_1,E_2)$ the
Banach space of all c.b.\ maps from $E_1$ to $E_2$.

We denote by $\Cal K$ the $C^*-$algebra 
 of all compact operators on $\ell_2$.

We will use repeatedly the notion of ``maximal'' operator space introduced
in [BP1], and further studied in [Pa6]. Let us recall its definition:\ let
$E$ be any normed space. Let $I$ be the class of all maps $u\colon \ B\to
B(H_u)$ with $\|u\|\le 1$ (and say $\dim H_u \le \text{card}(E)$). We let
$J\colon \ E\to \bigoplus\limits_{u\in I} B(H_u)$ be the isometric
embedding defined by $J(x) = \bigoplus\limits_{u\in I} u(x)$. Then,
$\max(E)$ is defined as the operator space $J(E) \subset
B\left(\bigoplus\limits_{u\in I}H_u\right)$, and any operator space which
is of this form (up to complete isometry) is called ``maximal''.

The ``maximal'' operator spaces are characterized by the property that,
for any linear map $u\colon \ E\to B(H)$ we have $\|u\|_{cb} = \|u\|$. The
following slightly more explicit description of their operator space
structure is often useful:\ for any $n$ and any $x$ in $M_n(\max(E))$ we
have $\|x\|<1$ iff, for some integer $N$, there is a diagonal matrix $D$
in $M_N(E)$ and scalar matrices $\beta\in M_{n,N}$ and $\gamma\in M_{N,n}$
such that 
$$x = \beta D\gamma\quad\text{\rm and}\quad  \|\beta\|\,
\|D\|\, \|\gamma\|<1.\tag0.4$$
We refer the reader to
[Pa6] for more information on this.

\bigskip

We now review the contents of this paper, section by section.

\n In section 1, we introduce the notion of ``similarity
setting'' which allows us to unify the various similarity
problems that we wish to consider. A similarity setting is
a triple $(i,E,\Cal A)$ where $E$ is an operator space,
$\Cal A\subset B(H)$ a unital subalgebra and $i\colon \
E\to \Cal A$ is an injective linear map with $\|i\|_{cb}
\le 1$, such that $\Cal A$ is generated by $i(E)$.

Given such a setting, for any $c\ge 1$ we construct the enveloping unital
operator algebra $\tilde A_c$ which contains $\Cal A$ as a dense
unital 
subalgebra and has the property that for any unital homomorphism $u\colon
\ \Cal A\to B(H)$, we have
$$\|ui\|_{cb(E,B(H))}\le c \Leftrightarrow \|u\|_{cb(\tilde A_c,B(H))}\le
1.$$
In particular, when $c=1$, $u$ is completely contractive on $\tilde A_1$
iff it is completely contractive when restricted to $E$. We also introduce in \S 
1,
the universal unital operator algebra (denoted by $OA(E)$) of an arbitrary
operator space $E$. The inclusion $E\to OA(E)$ can be viewed as the
``maximal'' setting involving $E$. The main result of \S 1 is Theorem~1.7
which gives an alternate description of $\tilde A_c$ as a canonical
quotient of $OA(E)$. This is the crucial tool used in \S 2,
where we present our theory of the similarity degree $d$ of
a setting $(i,E,\Cal A)$. This degree $d$ is defined as
the smallest number $\alpha\ge 0$ with the following
property:\ there is a constant $K$ such that for any $c\ge
1$ and any unital homomorphism $u\colon \ \Cal A\to
B(H)$ with $\|ui\|_{cb} \le c$ there is an invertible
$S\colon \ H\to H$ such that $e\to S^{-1}ui(e)S$ is
completely contractive and $\|S^{-1}\|\, \|S\| \le
Kc^\alpha$.

We prove (see Corollary 2.7) that $d$ is an integer and that the preceding
property still holds for $\alpha=d$.

In \S 3, we apply this to uniformly bounded group
 representations on a discrete group $G$, we denote
the degree in this case by $d(G)$, and we
prove the above Theorem~0.1, which implies that $d(G)\le 2$
iff $G$ is amenable.  We
actually prove a stronger version involving the space of
``coefficients'' of uniformly bounded (u.b.\ in short)
representations.

This is proved by applying \S 2 to the following setting:\ $E = \ell_1(G)$
with its (usual)  maximal  operator space structure 
(which also can be defined  by duality
with $c_0$,  \cf  [ER, BP1]), and $\Cal A\subset C^*(G)$
is the image of $\ell_1(G)$ under the canonical map from
$\ell_1(G)$ into $C^*(G)$.

 In \S 4, we come to the most natural ``setting'':\ we consider
a unital operator algebra $A\subset B(H)$ and we let $\Cal
A=A$ and $E = \max(A)$ in the sense of [BP1]. Then we
denote by $d(A)$ the corresponding degree. We give a
number of characterizations of this number.

In \S 5, we investigate the class of uniform algebras, i.e. $A\subset
C(K)$ ($K$ compact), $A$ is unital and separates the point of $K$. In
analogy with the group case, we prove that 
there is a constant $C$ such that
any unital homomorphism $u\: A\to B(H)$ satisfies
$\|u\|_{cb}\le C\|u\|^2$  iff $A=C(K)$, or
equivalently (by a result of Sheinberg [Sh]) iff $A$ is amenable. 

In \S 6, we turn to $C^*$-algebras and prove an analogous result (with the
assumption that $A$ has sufficiently many semi-finite representations):\
$d(A)\le 2$ iff $A$ is nuclear or equivalently (by results of Connes and
Haagerup, see [H4]) iff $A$ is amenable.

Finally, in \S 7, we give a slightly expanded 
version of some of
Blecher and Paulsen's results in [BP2]. We give, as an illustration,  an 
apparently new
characterization of the elements of the space $B(G)$
formed of the coefficients of unitary representations of
a discrete group $G$, to be compared with the case
 of uniformly bounded
representations treated in Theorem 1.12.

\head \S 1. Enveloping operator algebras. Preliminary results\endhead

\n It will be convenient to work in the following very general setting:\ we
give ourselves a unital algebra $\Cal A$ together with a linear subspace
$E\subset \Cal A$. We assume that $E$ is given with an operator space 
structure.
We will denote by
$$i\colon \ E\to \Cal A$$
the inclusion mapping. Moreover, we assume that the unital
algebra generated by $i(E)$ is the whole of $\Cal A$.

\n In addition,  we assume that $\Cal A$ can be faithfully
represented in $B(H)$ for some Hilbert space $H$ by a
unital 
representation $u_0\colon \ \Cal A\to B(H)$ such that
$\|u_0i\|_{cb}\le 1$. 

\n We will then say that the triple $(i,E,\Cal A)$ is a ``similarity setting".

Given such a setting, we can define for
any $c\ge 1$ the enveloping unital operator algebra
$\tilde A_c$ as follows:

\n Consider the family $\Cal C_c$ of all unital
homomorphisms $$u\colon \ \Cal A\to B(H_u)$$
with $H_u$ a Hilbert space, such that
$$\|ui\|_{cb}\le c.$$
Then we equip $\Cal A$ with the norm
$$\|a\|_c = \sup_{u\in \Cal C_c} \|u(a)\|.$$
Note that $\|a\|_c < \infty$ since $u$ is a homomorphism and $i(E)$
generates $\Cal A$. Moreover, since $u_0\in \Cal C_c$, we indeed have a 
norm.
We denote by $\tilde A_c$ the completion of $\Cal A$ for this norm. Clearly we
have an isometric unital homomorphism
$$\tilde A_c \subset \bigoplus_{u\in \Cal C_c} B(H_u)$$
which allows us to consider from now on $\tilde A_c$ as a unital operator
algebra (and a fortiori as an operator space).

Note that whenever $1\le c \le d$ we have $\Cal C_c \subset \Cal C_d$
hence we have a completely contractive unital homomorphism
$i_{c,d}\colon  \ \tilde A_d\to \tilde A_c$ with 
$$\|i_{c,d}\colon  \ \tilde A_d\to \tilde A_c
\|_{cb}\le 1.\tag1.0$$
 Note that $\tilde A_c$ is
characterized by the following property: \medskip
\roster
\item"{(1.1)}" any unital homomorphism $u\colon \ \Cal A\to B(H)$
such that $\|ui\|_{cb}\le c$ admits a unique  extension
$\tilde u\colon \ \tilde A_c \to B(H)$ with $\|\tilde
u\|_{cb(\tilde A_c, B(H))} \le 1$.
\endroster
\medskip

\n In this general setting, we wish to study the following\medskip

\n {\bf Similarity Property:}
\ For each $u$ in $\bigcup\limits_{c>1}\Cal C_c$, there is an invertible
operator $S\colon \ H\to H$ (= a similarity) such that the homomorphism
$$u_S\colon \ a\to S^{-1} u(a)S$$
satisfies $\|u_S i\|_{cb}\le 1$ (or equivalently is in
$\Cal C_1$).

\n  As we will see in the examples below, our
setting contains a number of fundamental similarity
problems:\ when $\Cal A$ is a group algebra (i.e.\
$\Cal A=\ell_1(G)$) or when $\Cal A$ is a $C^*$-algebra, or when $\Cal A$
is the disc algebra.\medskip

\example{Example 1.1} Let $G$ be a discrete group. Let $\Cal A$ be the group
algebra of $G$, i.e.\ $\Cal A = \ell_1(G)$ equipped with convolution. Let
$\Gamma\subset G$ be a set of generators for $G$ and let
$$E=\ell_1(\Gamma).$$
In this situation, it is easy to check that $\tilde A_1  = C^*(G)$ the
``full'' $C^*$-algebra of $G$ (= the enveloping $C^*$-algebra of
$\ell_1(G)$). Then the similarity property in this context means that for
any group representation $u\colon \ G\to B(H)$ such that
$\sup\limits_{\gamma\in \Gamma} \|u(\gamma)\|\le c$ there is a similarity
$S\colon \  H\to H$ such that $\sup\limits_{t\in G} \|S^{-1}u(t)S\|\le 1$. We 
study
this problem in section~3.\medskip
\endexample

\example{Example 1.2} Let $G={\nat}$ (a discrete semi-group can also be
discussed), let $E=\ell_1({\nat })$ and let $A(D)$ the
disc algebra with the natural contractive inclusion
$i\colon \ \ell_1({\nat })\to A(D)$.
We let $\Cal A=i(\ell_1({\nat }))$.
We equip
$\ell_1({\nat })$ with its maximal operator space
structure, so that for any map $v\colon \ \ell_1({\nat })\to B(H)$
we have $\|v\| = \|v\|_{cb}$.
 Consider a unital
homomorphism $u\colon \ A(D)\to B(H) $ such that $\|ui\| =
\|ui\|_{cb} \le c$, and let $T=u(z)$.  Then $T$ is a power
bounded operator and $$\|ui\| = \|ui\|_{cb} = \sup_{n\ge 1}
\|T^n\|.$$ Since there are power bounded operators which
are {\it not\/} similar to contractions ([Fo, Le]), the similarity
property does {\it not\/} hold in this case. 
\endexample

\example{Example 1.3} Let $\Cal A=A(D)$ the disc algebra and let $E=A(D)$
equipped
with its ``maximal'' operator space structure, $i$ being the identity on
$A(D)$. Then consider $u\colon \ \Cal A\to B(H)$ such that $\|ui\|_{cb} \le c$
and let $T= u(z)$. Here $ui$ is c.b.\ iff $T$ is
``polynomially bounded''. Moreover $\|ui\|_{cb}\le c$
holds iff we have $$\|P(T)\| \le c\|P\|_\infty$$
for any polynomial $P$.
\endexample

The similarity problem in this case is a well known  problem usually
attributed to Halmos.
The problem was solved by a counterexample in [P5].
 Analogous questions can
be formulated for any uniform algebra.
We will return to this topic in \S 5.\medskip

\example{Example 1.4} Let $\Cal A$ be a $C^*$-algebra
 and let $E=\max(\Cal A)$, with $i$
 again equal to the identity. Then
the similarity problem reduces  again to a well known open
problem raised by Kadison [Ka]:
\endexample

\n is every bounded unital
homomorphism $u\colon \ \Cal A\to B(H)$ similar to a
$*$-representation?

We discuss the $C^*$-algebra setting in \S 6. \vskip24pt 

Let $E$ be an arbitrary operator space. We wish to define the ``free unital
operator algebra'' associated to $E$. One way to define it is as follows.
We consider the free unital (noncommutative) algebra $\Cal P(E)$
associated to $E$ (equivalently, this is the tensor algebra
over $E$). The elements of $\Cal P(E)$ may be described as the
vector space of formal sums
$$P = \lambda_01 + \sum \lambda_{i_1} e^1_{i_1} +
\sum_{i_1i_2} \lambda_{i_1i_2} e^2_{i_1}e^2_{i_2}
 +\cdots+ \sum_{i_1\ldots i_N} \lambda^N_{i_1\ldots i_N} e^N_{i_1}
\ldots e^N_{i_N},\tag1.2$$
with $\lambda_0,\lambda_{i_1},\ldots, \lambda^N_{i_1\ldots i_N} \in {\comp}$ and 
with $e^1_i, e^2_{i_1},\ldots$ all in $E$, equipped with the
``free'' product operation.

Grouping terms, we may rewrite (1.2) as
$$P = P_0+P_1+\cdots+ P_N\tag1.3$$
with $P_0,P_1,\ldots, P_N$ ``homogeneous'', i.e.
$$P_N = \sum_{i_1\ldots i_N} \lambda^N_{i_1\ldots i_N}
e^N_{i_1}\ldots e^N_{i_N}.\tag1.3'$$
We will denote by $E^{(N)}$ the linear subspace of $\Cal P(E)$
spanned by all elements of the form (1.3)'. When $N=0$ we define
by convention
$$E^{(0)}=\comp 1.$$
The space $E^{(1)}$ is just $E$ viewed as a subset of $\Cal P(E)$. 
Then consider the family $J$ of all the mappings $v\colon \ E\to B(H_v)$
with $\|v\|_{cb}\le 1$. Let 
$$v(P) = \lambda_0I + \sum_i \lambda^1_iv(e^1_i) +\cdots+ \sum_{i_1\ldots
i_N} \lambda^N_{i_1i_2\ldots i_N} v(e^1_{i_1})\ldots
v(e^N_{i_N})$$ and
$$\|P\| = \sup_{v\in J} \|v(P)\|.$$
We will denote by $OA(E)$ the completion of $\Cal P(E)$ for this norm.
(The fact that it is a norm
 easily follows from (1.10) and (1.14) below.)
Clearly we have $\|PQ\| \leq \|P\|\, \|Q\|$ for all $P,Q$
in $\Cal P(E)$, hence we have a unital Banach algebra
structure on $OA(E)$.

\n We   denote by $OA_N(E)$ the closure in $OA(E)$ of
all the elements of the form  (1.3).

\n Moreover, we   denote by $E_N$ 
the closure in $OA(E)$ of the linear subspace $E^{(N)}$.

By construction, we have a natural embedding
$$OA(E) \subset \bigoplus_{v\in J} B(H_v)$$
which allows us to consider from now on $OA(E)$ as a
unital  operator algebra (and a fortiori as an operator
space) containing $E$ completely isometrically. This
operator space structure can be described as follows:\
consider an element $G$ in $\Cal K\otimes \Cal P(E)$.
Clearly $G$ can be written (for some $N$) as a finite sum
of the following form $$G = \lambda_0 \otimes 1 + \sum_i
\lambda^1_i\otimes e^1_i +\cdots + \sum_{i_1\ldots i_N}
\lambda^N_{i_1i_2\ldots i_N} \otimes e^N_{i_1} \ldots
e^N_{i_N}.\tag1.4$$
 with $\lambda^1_1,\ldots,
\lambda^N_{i_1\ldots i_N} \in \Cal K$ and
$$e^1_i,\ldots, e^j_{i_j}\ldots, e^N_{i_1},\ldots,
e^N_{i_N}\in E.$$ For short we will also write this as $$G
= G_0 +\cdots+G_N,\ \ G_j\in \Cal K\otimes E^{(j)} \tag1.5$$
 Then the
following formula encodes the operator space structure of $OA(E)$:
$$\align
\|G\|_{\Cal K\otimes_{\text{\rm min}} OA(E)} &=
\sup_{v\in J} \left\|\lambda_0 \otimes I + \sum_i
\lambda^1_i \otimes v(e^1_i)\right.\\
 &\left.\quad
+\cdots+ \sum_{i_1\ldots i_N} \lambda^N_{i_1\ldots i_N}
\otimes v(e^N_{i_1}) v(e^N_{i_2}) \ldots v(e^N_{i_N})
\right\|_{\Cal K \otimes_{\text{\rm min}} B(H_v)}.\tag1.6
\endalign$$
 This
algebra $OA(E)$ is characterized by the following (easily
verified) property:\medskip

\roster
\item"{(1.7)}" Let $B$ be any unital operator algebra. For any $v\colon \
E\to B$ with $\|v\|_{cb}\le 1$ there is a unique unital homomorphism $\hat
v\colon \ OA(E)\to B$ extending $v$ such that $\|\hat v\|_{cb} \le
1$.
\endroster
\medskip

\n (Note that actually the extension $\hat v$ is the restriction of a
$C^*$-algebra representation.)

For instance, we may consider $B = OA(E)$ and $v_z\colon \ E\to OA(E)$
defined by
$$ v_z(e)  = ze,\leqno \forall\ e\in E $$
where $z\in {\comp }$ with $|z|\le 1$. Then by construction of $OA(E)$,
$\|v_z\|_{cb} \le 1$ hence there is a unique unital homomorphism $\hat
v_z\colon \ OA(E)\to OA(E)$ extending $v_z$ and such that $\|\hat
v_z\|_{cb}\le 1$. We will use the notation
$$\omega(z) = \hat v_z.$$
Then it is easy to check (since $\hat v_z$ is a homomorphism extending
$v_z$) that if $P$ is as in (1.3) we have
$$\omega(z)P = P_0 +zP_1 +\cdots+ z^NP_N.$$
Similarly, if $G$ is as in (1.5) we have
$$(I_{\Cal K}\otimes \omega(z)) (G) = G_0 + zG_1 +\cdots+ z^NG_N.\tag1.8$$
It will be useful to record here the following fact.

\proclaim{Lemma 1.5} Each $P$ in $\Cal P(E)$ can be written in a unique
way as
$$P = P_0+P_1 +\cdots+ P_N,\tag1.9$$
for some integer $N$ with $P_j \in E^{(j)}$ for all $j\ge 0$. If we define
$$Q_j(P)  = P_j$$
(i.e.\ $Q_j(P) = 0$ $\forall j>N$) then $Q_j$ extends to a complete
contraction from $OA(E)$ onto $E^{(j)}$.
\endproclaim

\demo{Proof} By the preceding formula (1.8) we have
$$G_j = \int \bar z^j[I_{\Cal K} \otimes\omega(z)](G) dm(z)$$
where $m$ denotes normalized Haar measure on $\{z\mid |z|=1\}$. By
convexity this yields
$$\|G_j\|_{\text{\rm min}} \le \|G\|_{\text{\rm min}}.\tag1.10$$
Since $G_j = (I_{\Cal K} \otimes Q_j)(G)$, this shows that $\|Q_j\|_{cb}\le
1$, and it also shows the unicity of the expression (1.9).\qed
\enddemo

\remark{Remark} Another description of $OA(E)$ is as follows:\ we consider the
$C^*$-algebra $C^*<E>$ constructed in [Pe]. This is characterized by the
property that for  any $v\colon \ E\to B$ ($B$ any $C^*$-algebra) with
$\|v\|_{cb}\le 1$ there is a representation
 $\pi\colon \ C^*<E>\to B$
extending $v$. Then we can define $OA(E)$ as the unital (non-selfadjoint)
operator algebra generated by the elements of $E$ in
the unitization of  $C^*<E>$.
\endremark

We now introduce the ``product map" $\pi_1$ and a whole family
of deformations $\pi_z$. 

\n Consider $z\in \comp$ with $|z|\le 1$ and let $c=1/|z|$.
We can define a unital homomorphism
$$\pi_z\colon \ OA(E) \to \tilde A_c$$
as follows:

\n Let $V_z = z i\colon \ E\to \tilde A_c$. Then $\|V_z\|_{cb}\le 1$.
Therefore, by (1.7) (since $\tilde A_c$ is an operator algebra) there is a 
unique
unital homomorphism
$$\pi_z\colon \ OA(E)\to \tilde A_c$$
extending $V_z$ and such that $\|\pi_z\|_{cb}\le 1$.

We will need the following simple observation.
\def\ie{ {\it i.e.}\ }

\proclaim{Lemma 1.6} Consider our usual similarity setting $(i,E,A)$.
 Assume that $E$ contains
the unit element $1_{\Cal A}$ of $\Cal A$. Let $0\le j<N$. Then for any $g$ in 
$\Cal K\otimes E^{(j)}$, there is $g'$ in $\Cal K\otimes E^{(N)}$ such
that $$(I_{\Cal K} \otimes \pi_1)(g-g') = 0\tag1.11$$
$$\|g'\|_{\text{\rm min}} \le \|1_{\Cal A}\|^{N-j}_E \|g\|_{\text{\rm min}}. 
\tag1.12$$
\endproclaim

\demo{Proof}   We introduce the map $V\colon \ OA(E)\to OA(E)$ which is simply 
the
left multiplication by $1_{\Cal A} $, \ie $V(x)=x\otimes 1_{\Cal A}$.    Clearly 
we have
$$ V[E^{(j)} ]\subset E^{(j+1)} \leqno \forall j\ge 0$$
hence
$$V^{N-j}[E^{(j)}] \subset E^{(N)}.$$
Moreover for any $P$ in $\Cal P(E)$ we have clearly $\pi_1(V(P)) =
\pi_1(P)$ hence
$$\pi_1(V^{N-j}(P)) = \pi_1(P).$$
Similarly for any $g$ in $\Cal K\otimes E^{(j)}$, let
$$g' = (I_{\Cal K} \otimes V^{N-j})(g).$$
Then we clearly have (1.11) and (1.12). (Indeed, since $OA(E)$ is an operator
algebra \break $V\colon \ x\to 1_{\Cal Ax}$ has $cb$ norm $\le
\|1_{\Cal A}\|_E$.) \qed
\enddemo

 We now come to our
first result.

\proclaim{Theorem 1.7} Let $c\ge 1$ and $z=1/c$. The 
mapping $\pi_z$ is a completely contractive 
surjection from $OA(E)$ onto $\tilde A_c$. Moreover, it induces canonically
a completely isometric isomorphism
$$\sigma_z\colon \ OA(E)/\ker(\pi_z) \to \tilde A_c,$$
so that, if $Q_z\colon \ OA(E)\to OA(E)/\ker(\pi_z)$ denotes the
canonical surjection, we have
$$\pi_z = \sigma_z Q_z.$$
More precisely, for any $f$ in $\Cal K\otimes \Cal A$ with $\|f\|_{{\Cal K}
\otimes_{\text{\rm min}} \tilde A_c}<1$, there is $F$ in
 $\Cal K
\otimes \Cal P(E)$ with $\|F\|_{\Cal K\otimes_{\text{\rm min}}
OA(E)} < 1$ such that $$I_{\Cal K} \otimes \pi_z(F) =f.$$
\endproclaim

\demo{Proof} Note that $\pi_z$ is characterized as the unital homomorphism such 
that
$$\pi_z(e) = zi(e),\leqno \forall \ e\in E$$
if we view $E\subset OA(E)$ and $\Cal A\subset \tilde A_c$. By construction we
have $\|\pi_z\|_{cb}\le 1$. On the other hand, note that $\pi_z(OA(E))$
contains $i(E)$ and is a subalgebra, hence it contains $\Cal A\subset \tilde
A_c$ since we assume that $i(E)$ is generating. Therefore we can define a
mapping
$$u_z\colon \ \Cal A\to OA(E)/\ker(\pi_z)$$
simply by setting
$$  u_z(a) = (\sigma_z)^{-1}(a).\leqno \forall\ a \in \Cal A\subset
\tilde
A_c$$
This is clearly a unital homomorphism. Moreover we have $u_zi(e) =
Q_z(z^{-1}e)$ since $\sigma_z Q_z(z^{-1}e) =
\pi_z(z^{-1}e) = i(e)$. Since $u_zi = z^{-1}Q_{z|E}$ it
follows that $$\|u_zi\|_{cb}\le |z|^{-1} =c.$$
Hence by the defining property of $\tilde A_c$, (since
$OA(E)/\ker(\pi_z)$ is an operator algebra [BRS])
there is a unique unital homomorphism $\tilde u_z\colon \
\tilde A_c \to OA(E)/\ker(\pi_z)$ such that
 $\|\tilde u_z\|_{cb}\le 1$. Moreover $\tilde u_{z|\Cal A} = u_z$ hence we
have $\|(\sigma_z)^{-1}(a)\| = \|u_z(a)\| \le
\|a\|_{\tilde A_c}$ for any $a$ in $\Cal A$. By the density of
$\Cal A$ in $\tilde A_c$, it follows that $\sigma_z$ is a
surjective isometry, and $\tilde u_z = (\sigma_z)^{-1}$,
so that finally $$\|(\sigma_z)^{-1}\|_{cb} = \|\tilde
u_z\|_{cb} \le 1.$$ Thus $\sigma_z $ is a complete
isometry.

\n For the last assertion in Theorem~1.7, we need an obvious extension of the
main result in [BRS] to non-closed unital subalgebras of $B(H)$, as
follows:\ let $\Cal P$ be a
 unital subalgebra of $B(H)$ and let $I\subset
\Cal P$ be a 2-sided ideal which is closed in $\Cal P$. Then the quotient 
space $Z = \Cal P/I$ can
be equipped with a (noncomplete)
 operator space structure by Ruan's theorem
[R]. Consequently, its completion $\widehat Z$ can be equipped with a
(complete this time) operator space structure.

\n On the other hand, $\widehat
Z$ is clearly a unital Banach algebra, and it is easy to check that the
product mapping $\widehat Z \otimes_{h} \widehat Z\to \widehat Z$ is a
complete contraction. Hence by [BRS] there is a completely
isometric unital homomorphism from $\widehat Z$ into
$B(H)$ for some Hilbert space $H$.

\n Returning to the situation in Theorem~1.7, let $\Cal P = \Cal P(E)$
(equipped with the operator space structure induced by $OA(E)$) and $I =
\Cal P(E) \cap \ker(\pi_z)$. Let us denote $Z_c = \Cal P/I$ in
this case. Clearly, the restriction of $\pi_z$ to $\Cal P(E)$ induces a
completely contractive unital homomorphism $\hat\sigma_z\colon \ \widehat
Z_c \to \tilde A_c$, which is injective on $Z_c$.

\n Since we assume that $i(E)$ generates $\Cal A$, we have $\hat\sigma_z(Z_c) 
=
\pi_z(\Cal P(E))=\Cal A$. Thus the inverse of $\hat\sigma_{z|Z_c}$ defines a
homomorphism $u_z\colon \ \Cal A\to \widehat Z_c$ such that $\|u_zi\|_{cb}\le
c$, and repeating the preceding argument we obtain that $\hat \sigma_z$
must be a complete isometry from $\widehat Z_c$ onto $\tilde A_c$.\qed
\enddemo

The next result is a simple reformulation of Paulsen's results in [Pa4].

\proclaim{Proposition 1.8} Let $K\ge 0$ be a constant. Let $E$ be an
operator space, $\Cal A$ a unital algebra and $i\colon \ E\to \Cal A$ an 
injection
with $i(E)$ generating $\Cal A$. The following properties of a unital
homomorphism $u\colon \ \Cal A\mapsto B(H)$ are equivalent, for each fixed 
$c\ge
1$.\medskip
\roster
\item"{(i)}" $u$ extends to a c.b. homomorphism
$$\tilde u\colon \ \tilde A_c\to B(H) \quad \text{with}\quad \|\tilde
u\|_{cb}\le K.$$
\item"{(ii)}" There is an invertible operator $S\colon \ H\to H$ with
$\|S\|\, \|S^{-1}\| \le K$ such that the map $u_S\colon \ \Cal A\to B(H)$
defined by $u_S(a) = S^{-1}u(a)S$ extends completely contractively to
$\tilde A_c$.
\item"{(iii)}" There is an invertible operator $S$ with $\|S\|\, \|S^{-1}\|
\le K$ such that
$$\|u_S i\|_{cb}\le c.$$
\endroster
\endproclaim

\demo{Proof} The equivalence (i) $\Leftrightarrow$ (ii) is exactly Paulsen's 
result
[Pa4]. Clearly (ii)~$\Rightarrow$~(iii) holds since by construction
$\|i\colon \ E\to \tilde A_c\|_{cb} \le c$. Now assume (iii). By the defining
property of $\tilde A_c$, $u_S\colon \ \Cal A\to B(H)$ admits an 
extension
$\tilde u_S\colon \ \tilde A_c \to B(H)$ with $\|\tilde u_S\|_{cb}\le 1$.
 Hence we obtain (ii).\qed
 \enddemo

\remark{Remark 1.9} Let $E$ be any operator space. Consider the iterated
Haagerup tensor product $X = E\otimes _h E\ldots \otimes_h E$ ($N$-times).
Let $x$ be arbitrary in $\Cal K\otimes E\otimes\cdots\otimes E$. Then $x$
can be written as a finite sum
$$x = \sum_{i_1i_2\ldots i_N} \lambda_{i_1i_2\ldots i_N} \otimes x^1_{i_1}
\otimes\cdots\otimes x^N_{i_N}$$
with $\lambda_{i_1\ldots i_N} \in \Cal K$ and $x^1_{i_1},\ldots,
x^N_{i_N} \in E$. It is proved in [CES] that we have
$$\|x\|_{\Cal K\otimes_{\text{\rm min}}X} = \sup\left\{\left\| \sum_{i_1\ldots
i_N} \lambda_{i_1i_2\ldots i_N} \otimes \sigma^1(x^1_{i_1})
\sigma^2(x^2_{i_2}) \ldots \sigma^N(x^N_{i_N})\right\|_{\text{\rm 
min}}\right\}$$
where the supremum runs over all possible choices of $H$ and of complete
contractions $\sigma^1\colon \ E\to B(H),\ldots, \sigma^N\colon \ E\to
B(H)$. We claim that actually this supremum is attained when
$\sigma^1,\ldots, \sigma^N$ are all the same, more precisely we have
$$\|x\|_{\Cal K\otimes_{\text{\rm min}}X} = \sup\left\{\left\|
\sum_{i_1\ldots i_N} \lambda_{i_1i_2\ldots i_N} \otimes
\sigma(x^1_{i_1}) \sigma(x^2_{i_2}) \ldots
 \sigma(x^N_{i_N})\right\|_{\text{\rm min}}\right\}
\tag1.13$$ 
where the supremum runs over all possible
$H$ and all complete contractions $\sigma\colon \ E\to
B(H)$.\medskip
\endremark

\n Indeed, this follows from a trick already used by Blecher in [B1] and which 
seems to originate in Varapoulos's
paper [V]. The trick consists in replacing $\sigma^1,\ldots, \sigma^N$ by
the single map $\sigma\colon \ E\to B(\underbrace{H\oplus\cdots\oplus
H}_{N+1\ \text{\rm times}})$ of the form
$$\sigma(e) = \left(\matrix 0&&\sigma^1(e)\\
&&\ddots&&\bigcirc\\
&&&\ddots\\
&&&&\sigma^N(e)\\
\bigcirc\cr &&&&0\endmatrix\right).$$
(More precisely $\sigma(e)$ is the $(N+1)\times (N+1)$ matrix having
$(\sigma^1(e),\ldots, \sigma^N(e))$ above the main diagonal and zero
elsewhere).

\n Then it is easy to check that $\|\sigma\|_{cb} =
\sup_j\|\sigma^j\|_{cb}$ and $\forall x^1,\ldots, x^N\in E$
$$[\sigma(x^1)\ldots \sigma(x^N)]_{1,N+1}  =
\sigma^1(x^1)\sigma^2(x^2) \ldots \sigma^N(x^N).$$
>From this our claim immediately follows. Note that our claim shows for
instance that in the case $E  = \max(\ell^n_1)$, the space $E\otimes_h
\cdots \otimes_h E$ ($N$ times) can be identified completely isometrically
with a subspace of $C^*({\bold F}_n)$ (here ${\bold F}_n$ is the free group with 
$n$
generators). Namely the subspace spanned by all products
$$\{U_{i_1}U_{i_2}\ldots U_{i_N}\mid 1\le i _1 \le n,\ldots, 1\le i_N \le
n\},$$
where $U_1,\ldots, U_n$ denote the free
unitary generators of $C^*({\bold F}_n)$. Although this
useful  fact might
have been observed by others,
it does not seem to have been recorded into print.

\proclaim{Proposition 1.10} Let $E$ be any operator space. Consider $E$ as
embedded into $OA(E)$. Fix $N\ge 1$, recall that we
 denote by $E_N$ the closed
subspace of $OA(E)$ 
generated by all products of the form $x_1\cdot x_2\ldots
x_N$ with $x_i\in E$. Then the natural ``product'' mapping
$$T_N\colon \ E\otimes_h\ldots \otimes_h E\to  {E_N}$$
which takes $x_1\otimes \cdots\otimes x_N$ to $x_1\cdot x_2\ldots \cdot
x_N$ is a completely isometric isomorphism.
\endproclaim

\demo{Proof} For simplicity let us denote $X =E\otimes_h\cdots\otimes_h E$ ($N$
times). Since the algebraic tensor product
$E\otimes\cdots \otimes E$ is dense in $X$ and similarly for ${E_N}$, it
suffices to prove that for any element $G$ in $\Cal K \otimes E \otimes
\cdots \otimes E$ we have
$$\|G\|_{\Cal K\otimes _{\text{\rm min}}X} = \|(I_{\Cal K}\otimes T_N)
(G)\|_{\Cal K\otimes_{\text{\rm min}}{E_N}}.\tag1.14$$
But this is immediate by (1.13) and (1.6).\qed
\enddemo

We now record here several consequences of Theorem~1.7 
(inspired by Peller's results for the category of
$Q$-algebras in [Pel, Prop. 4.2 and 4.3]).

\proclaim{Corollary 1.11} Let $W\colon\  \tilde A_c \to B(H)$ be a
linear mapping. Let $z=1/c$. Then
$$\|W\|_{cb(\tilde A_c, B(H))} = \|W\pi_z\|
_{cb(OA(E), B(H))}.\tag1.15$$
Moreover, for any linear map $w\colon \ \Cal A\to B(H)$, the following
assertions are equivalent:\medskip
\roster
\item"{\text{\rm (i)}}" For some $c\ge 1$ and some $K\ge 0$, $w$ extends to a 
$c.b.$ 
map
$\widetilde w\colon \ \tilde A_c \to B(H)$ with $$\|\widetilde w\|_{cb} \le
K.$$ 
\item"{\text{\rm (ii)}}" There are constants $c'\ge 1$ and $K'\ge 0$ such that, 
for any
$N\ge 1$, the mapping $w_N\colon \ E\otimes_h\cdots \otimes_h E\to B(H)$
defined by $w_N(x_1 \otimes\cdots\otimes x_N) = w(x_1  \ldots  x_N)$
satisfies
$$\|w_N\|_{cb(E\otimes_h\cdots\otimes_h E, B(H))} \le K'(c')^N.$$
\item"{\text{\rm (ii)'}}" There are constants $c'\ge 1$ and $K'\ge 0$ such that, 
for any
$N\ge 1$,
there are bounded linear mappings $u_i\colon \ E\to B(H_{i+1},
H_i)$, with  $\|{u_i}\|_{cb}\le c'$ for all $i$, 
where $H_i$ are Hilbert spaces with $H_{N+1}=H$ and $H_1=H$, such
that
$$w(x_1x_2\ldots x_N) = K' u_1(x_1) u_2(x_2) \ldots
u_N(x_N).\leqno \forall\ x_1,\ldots, x_N\in E$$
\endroster
\endproclaim

\demo{Proof} The first part is an obvious consequence of the first assertion in
Theorem~1.7.
 We now prove the second part.
 The equivalence between (ii) and (ii)', 
with the same constants $K',\ c'$, is a particular case of
the well known factorization 
of completely bounded multilinear forms (\cf  [CS1-2, PaS]).
We now turn to the remaining equivalence. Assume (i). Then $\|\widetilde w
\pi_z\|_{cb}\le K$, hence $\|\widetilde w\pi_{z|E_N}\|_{cb}\le K$, but for
$x_1,\ldots, x_N$ in $E$ we have
$$\widetilde w\pi_z(x_1,\ldots, x_N) = c^{-N} w(x_1 \ldots  x_N) =
c^{-N} w_N(x_1\otimes \cdots\otimes x_N).$$
Hence by Proposition 1.10, we have
$$\|w_N\|_{cb} = c^N \|\widetilde w\pi_{z|E_N}\|_{cb}\tag1.16$$
hence
$$\|w_N\|_{cb} \le Kc^N.$$
This proves (i) $\Rightarrow$ (ii).

\n Conversely, assume (ii). Let $c>c'$ and $z=1/c$ as before. Then we have
by (1.16)
$$\|\widetilde w\pi_{z|E_N}\|_{cb}\le K'\left({c'\over c}\right)^N.
\leqno \forall\ N\ge 1$$
By Lemma 1.5, this implies using (1.15),
$$\|\widetilde w\|_{cb(\tilde A_c, B(H))} =
\|\widetilde w\pi_z\|_{cb} \le 1 + \sum_{N\ge 1} \|\widetilde
w\pi_{z|E_N}\|_{cb} \le 1+K' \sum_{N\ge 1} (c'/c)^N <\infty,$$
whence (i) with $K = 1+K' \sum\limits_{N\ge 1} (c'/c)^N$.\qed
\enddemo

We now illustrate the meaning of Corollary 1.11 in the group case.

\proclaim{Theorem 1.12} Let $G$ be a group (or merely a semi-group). Consider
a function $f\colon \ G\to B(H)$. The following assertions are
equivalent.\medskip
\roster
\item"{\text{\rm (i)}}" There is a uniformly bounded representation $\pi\colon \ 
G\to
B(H_\pi)$ and bounded operators $\xi\colon \ H_\pi\to H$ and $\eta\colon \
H\to H_\pi$ such that
$$f(t) = \xi\pi(t)\eta.\leqno \forall\ f\in G$$
\item"{\text{\rm (ii)}}" There are constants $K'\ge 0$ and $c'\ge 1$ 
such that, for each
$N\ge 1$, the function $f_N\colon\ G^N\to B(H)$ defined by
$$f_N(t_1,t_2,...,t_N)=f(t_1t_2\ldots t_N)$$
defines (with the obvious identification) an element
of $cb(\ell_1(G)\otimes_h\ldots\otimes_h\ell_1(G), B(H))$ 
($N$-fold tensor product) with norm $\le K'c'{}^N.$
\item"{\text{\rm (ii)'}}" There are constants $K'\ge 0$
 and $c'\ge 1$ such that, for each
$N\ge 1$, there are bounded mappings $\varphi_i\colon \ G\to B(H_{i+1},
H_i)$, with $\sup_{t\in G} \|\varphi_i(t)\|\le c'$ for all $i$, 
where $H_i$ are Hilbert spaces with $H_{N+1}=H$ and $H_1=H$, such
that
$$f(t_1t_2\ldots t_N) = K'\varphi_1(t_1) \varphi_2(t_2) \ldots
\varphi_N(t_N).\leqno \forall\ t_1,\ldots, t_N\in G$$
\endroster
\endproclaim

\demo{Proof} We merely apply Corollary 1.11 and Proposition~1.8 with $\Cal A=E =
\ell_1(G)$.     Note that, using the factorization of cb maps, it is easy to 
verify
that (i) holds iff the mapping
$t\to f(t)$ extends linearly to a mapping $\widetilde w\colon \ \tilde A_c \to 
B(H)$
 with $\|\widetilde w\|_{cb} \le
K.$ We leave the details to the reader.
\qed
\enddemo

\remark{Remark} Note that if (i) holds in the preceding statement with
$|\pi|\le c$ then we obtain (ii) with $K = \|\xi\|\, \|\eta\|$ and the
same number $c$. However, if (ii) holds we only obtain (i) with a
representation $\pi$ such that $|\pi|\le (1+\vp)c'$ (with $\vp>0$) and with 
$\|\xi\|\,
\|\eta\|\le K_\vp = 1+K \sum\limits_{N\ge 1} (1+\vp)^{-N}$. Indeed, these
are the constants appearing in the proof of Corollary~1.11. Nevertheless,
we will see below (see Corollary 7.8) that, in the particular case $c'=1$, we 
can get
rid of this extra factor $(1+\vp)$.
\endremark

\head\S 2. Main results\endhead

\n Let $E,\Cal A$ and $i\colon \ E\to \Cal A$ be our general setting as 
described in the
beginning of \S 1. We will assume that the following holds:

Every unital homomorphism $u\colon \ \Cal A\to B(H)$ such that $\|ui\|_{cb} <
\infty$ is similar to a homomorphism such that $\|ui\|_{cb}\le 1$, i.e.\
there is an invertible operator $S\colon \ H\to H$ such that the map $e\to
S^{-1}ui(e)S$ is completely contractive.

When this holds we will say that in this setting the similarity property
holds. We will need to carefully keep track of the constants involved in
this phenomenon.

\proclaim{Lemma 2.1} If the similarity property holds
then there is, for
each $c\ge 1$, a number $\Phi(c)$ such that,
 for any unital homomorphism
$u\colon \ \Cal A\to B(H)$  with $\|ui\|_{cb}\le c$, there is a similarity
$S\colon \ H\to H$ with $\|S\|\, \|S^{-1}\|\le \Phi(c)$ such that $e\to
S^{-1} ui(e)S$ is a completely contractive map from $E$ to $B(H)$.
\endproclaim

\demo{Proof} This is elementary. Just consider the unital
homomorphism  $U = \bigoplus\limits_{u\in \Cal C_c}
u$ and a similarity $\Cal S$ such that $\Cal S^{-1}U
\Cal S$ is contractive then restrict to the invariant
subspaces associated to each $u$ in $\Cal C_c$. We get
the announced bound with $\Phi(c) = \|\Cal S\|\, \|\Cal
S^{-1}\|$.\qed
\enddemo

The preceding lemma allows us to define the following parameter associated
to the similarity property
$$\Phi(i,c) = \sup_{u\in \Cal C_c} \inf\{\|S\|\, \|S^{-1}\|\}\tag2.1$$
where the infimum runs over all $S\colon \ H_u\to H_u$ invertible such that
$\|u_S i\|_{cb}\le 1$ where $u_S(a) = S^{-1} u(a)S$. When
the supremum is infinite, we write $\Phi(i,c) = \infty$ by
convention. Equivalently by Proposition~1.8 we have
$$\Phi(i,c) = \sup_{u\in \Cal C_c} \|u\|_{cb(\tilde A_1, B(H_u))}, \tag2.2$$
where by $\|u\|_{cb(\tilde A_1, B(H_u))}$ we mean that we compute the $cb$
norm of $u\colon \ \Cal A\to B(H_u)$ using the operator space
structure  induced by $\tilde A_1$ on $\Cal A$. Since $\Cal A$ is
dense in $\tilde A_1$, there is no risk of confusion.

By the definition of $\tilde A_c$ and by (2.2), when the similarity
property holds, then the natural map $\tilde A_c\to \tilde
A_1$ (which is always a complete contraction by (1.0)) is a
complete isomorphism and (2.2) can be rewritten as
$$\Phi(i,c) = \|\tilde A_1\to \tilde A_c\|_{cb}.\tag2.2'$$
 It will be convenient to introduce the following
notation $$\text{Sim}(u) = \inf\{\|S^{-1}\|\,
\|S\|\}\tag2.3$$
 over all $S$ such that $\|u_S
i\|_{cb}\le 1$. By convention we set $\text{Sim}(u)=+\infty$ if there is
no such $S$.

\proclaim{Lemma 2.2} Let $c>1$. Assume $\Phi(i,c) <\infty$. Then for any
$0<\theta < 1$ we have
$$\Phi(i,c) \le \Phi(i,c^\theta)^{1/\theta}.$$
\endproclaim

\demo{Proof} Consider $u$ in $\Cal C_c$ with $\text{Sim}(u)$ finite. By 
definition of 
$\Phi(i,c)$ for any
$\varepsilon>0$ $\exists S\colon \ H\to H$ invertible such that
$\|S^{-1}\|\ \|S\|\le \text{Sim}(u) +\varepsilon$ and $\|u_Si\|_{cb}\le 1$
where we again denote $u_S(a) = S^{-1} u(a)S$. Clearly we
can assume that $S$ is hermitian. We then invoke the three
lines lemma. Consider $z\in {\comp }$ with $0\le
\text{Re}(z) \le 1$ and $e\in E$.

\n If $\text{Re } z = 1$ we have $\|S^{-z}ui(e)S^z\|\le
\|e\|$, and if $\text{Re } z = 0$ we have $\|S^{-z}
ui(e)S^z\|\le c\|e\|$. Hence by the subharmonicity of
$$z\to \text{Log}\|S^{-z}ui(e)S^z\|$$
we have
$$\|S^{\theta-1}ui(e)S^{1-\theta}\|\le c^\theta\|e\|.$$
More generally, the same reasoning exactly yields that the map $v\colon \
\Cal A\to B(H)$ defined by $v(a) = S^{\theta-1}u(a)S^{1-\theta}$ satisfies
$\|vi\|_{cb}\le c^\theta$. In other words we have $v\in \Cal
C_{c^\theta}$ so that by definition of $\Phi(i,c^\theta)$ there is a
similarity $T\colon \ H\to H$ such that $\|T^{-1}\|\, \|T\|\le
\Phi(i,c^\theta)+\varepsilon$, and $\|v_Ti\|_{cb}\le 1$. This last
inequality implies that $e\to T^{-1}S^{\theta-1}ui(e)S^{1-\theta}T$ is a
complete contraction. Hence
$$\text{Sim}(u) \le \|T^{-1}S^{\theta-1}\|\, \|S^{1-\theta}T\| \le
\|T^{-1}\|\, \|T\| (\|S\|\, \|S^{-1}\|)^{1-\theta}$$
i.e.
$$\text{Sim}(u) \le (\Phi(i,c^\theta)+\varepsilon) (\text{Sim}(u) +
\varepsilon)^{1-\theta}.$$
Since this holds for all $u$ in $\Cal C_c$, we have
$$\Phi(i,c)\le (\Phi(i,c^\theta) + \varepsilon) (\Phi(i,c) +
\varepsilon)^{1-\theta}$$
now since $\varepsilon>0$ is arbitrary and $\Phi(i,c)<\infty$, we obtain
after division by $\Phi(i,c)^{1-\theta}$
$$\Phi(i,c)^\theta \le\Phi(i,c^\theta).$$
\qed
\enddemo

\remark{Remark} In case $\Phi(i,c)=\infty$ (for some $c>1$), the preceding
statement remains valid if there exist unital homomorphisms $u\colon\ \Cal
A\to B(H)$ with $\|ui\|_{cb}\le c$ for which $\text{Sim}(u)$ takes 
arbitrarily large {\it finite} values. (By
Haagerup's result in [H1] on cyclic homomorphisms, this condition is always
satisfied in the $C^*$-setting considered in \S 6.)
Then the preceding proof
shows that $\Phi(i,c')=\infty$ for {\it all} $c'>1$.
\endremark

\proclaim{Lemma 2.3} Let $\varphi\colon \  [1,\infty[ \to {\reel }_+$ be a
non-decreasing function such that $\varphi(c) \le
\varphi(c^\theta)^{1/\theta}$ whenever $0<\theta < 1$ and $c\ge 1$. Fix
$\alpha>0$. Then, if $\varphi(c)< c^\alpha$ for some $c>1$ we have
$\varphi(x) < x^\alpha$ for all $x\ge c$ and there is a constant $K$
such that $\varphi(x) \le K x^\alpha$ for all $x\ge 1$.
\endproclaim

\demo{Proof} Assume $c>1$ and $\varphi(c) = c^\beta$ with $\beta<\alpha$. Then 
let
$t
= 1/\theta$ with $1\le t<\infty$. Since $\varphi(x^t) \le \varphi(x)^t$ for
any $x\ge 1$ we have $\varphi(c^{t }) \le \varphi(c)^{t} = c^{\beta
t^n}$. Hence we have $\varphi(x) \le x^\beta$ whenever 
$x=c^{t}$ with
  $t\ge 1$. Equivalently this holds whenever $x\ge c$.
 The rest is obvious by monotonicity.\qed
\enddemo

Hence we have

\proclaim{Corollary 2.4} If the similarity property holds (i.e.\ if
$\Phi(i,c)<\infty$ for all $c\ge 1$) 
then there is a constant $K$ and an
exponent $\alpha\ge 0$ such that
$$ \Phi(i,c) \le Kc^\alpha.\leqno \forall c\ge 1$$
\endproclaim

We come now to our main result.

\proclaim{Theorem 2.5} Consider our usual setting $E,\Cal A,i$ and assume that
the similarity property holds. More precisely we assume that for some
constants $K>0$ and $\alpha>0$ we have $\Phi(i,c) \le Kc^\alpha$ for all
$c$ large enough.
Let $N$ be an integer with $\alpha < N+1$.
 Then the restriction of $\pi_1$
to $OA_N(E)$ is a complete surjection 
of $OA_N(E)$ onto $\tilde A_1$, i.e.\
there is a constant $K'$ such that, for any $f$ in $\Cal K \otimes_{\text{\rm
min}} \tilde A_1$ with $\|f\|_{\text{\rm min}} <1$, there is $\hat f$ in $\Cal K
\otimes_{\text{\rm min}} OA_N(E)$ with $\|\hat f\|_{\text{\rm min}} < K'$ such 
that
$(I_{\Cal K} \otimes \pi_1)(\hat f) = f$. Furthermore, this last property
implies that for some constant $K_1$ we have
$$\forall\ c\ge 1\qquad  \Phi(i,c) \le K_1 c^N.\tag2.4 $$
Finally, if $1_{\Cal A}$ belongs to $E$, then the   restriction of $\pi_1$
to $ {E_N}$ is a complete surjection of  $ {E_N}$ onto $\tilde A_1$.
\endproclaim

\demo{Proof of Theorem 2.5} Let us denote for simplicity
$$\Cal X_N = OA_N(E).$$
We first fix $c>1$ chosen large enough so
that
$$\sum_{j>N} Kc^{\alpha-j} < 1/2.\tag2.5$$
Note that since $N+1>\alpha$, this choice is possible. By the standard
iteration argument used in the proof of the open mapping theorem, it
suffices to prove the following.
\enddemo

\n {\bf Claim.} There is a constant $K''$ such that for
any $f$ in the open unit ball of $\Cal K \otimes_{\text{\rm
min}} \tilde A_1$, with $f\in \Cal K \otimes \Cal A$, there is
an element $\hat f$ in $\Cal K\otimes_{\text{\rm min}} \Cal
X_N$ with $\|\hat f\|_{\text{\rm min}} < K''$ and such that
$$\|(I_{\Cal K} \otimes \pi_1) (\hat f) -f\|_{\text{\rm min}} <
1/2.\tag2.6$$ 
By our assumption, we have a natural isomorphism
$$\varphi_z\colon \ \tilde A_1 \to \tilde A_c$$ 
which is the identity on $\Cal A$, with
$\|\varphi_z\|_{cb}\le Kc^\alpha$, $z=1/c$. Let $f$ be as
in our present claim. Then we have $$\|(I_{\Cal K} \otimes
\varphi_z)(f)\|_{\Cal K\otimes_{\text{\rm min}} \tilde A_c} <
Kc^\alpha.$$ Hence by Theorem 1.7 there is $g$ in
 $\Cal K \otimes_{\text{\rm min}} OA(E)$ such that
$$\|g\|_{\text{\rm min}} < Kc^\alpha \quad \text{and}\quad
 (I_{\Cal K} \otimes \pi_z)(g) =(I_{\Cal K} \otimes\varphi_z)(f).\tag2.7$$
Note that since $\varphi_z$ is the identity on $\Cal A$ and
$f\in \Cal K\otimes \Cal A$,
we may write $(I_{\Cal K} \otimes\varphi_z)(f)=f$.
  
\n We can assume that, for some $m$,
  $g$ is of the form $g=g_0 +\cdots+ g_m$
with $g_j \in \Cal K\otimes E^{(j)}$.
  By Lemma 1.5 we have
$$\forall \ j\le m \qquad\|g_j\|_{\text{\rm min} }\le \|g\|_{\text{\rm min}} < 
Kc^\alpha.
\tag2.8$$ Now
 let $G_j = (I_{\Cal K}\otimes
\pi_1)(g_j)\in \Cal K\otimes_{\text{\rm
min}\}tilde A_1}$. Note that
 $$f = (I_{\Cal K} \otimes \pi_z)(g)= 
\sum\nolimits^m_0
z^jG_j,$$
and $\|G_j\|_{\min}\le \|g_j\|_{\min} < Kc^\alpha.$
 Hence we have
$$
\left\|f-\sum\nolimits^N_0 z^jG_j\right\|_{\text{\rm min}}	\le
 \sum_{j>N} |z|^j \|G_j\|_{\text{\rm min}}
\le  \sum_{j>N} Kc^{\alpha-j},$$
therefore by (2.5),
 we obtain
$$\left\|f-\sum\nolimits^N_0 z^jG_j\right\|_{\text{\rm min}} < 1/2.$$
 Let $\hat f=\sum\nolimits^N_0 z^j g_j$, then by (2.8)
we have $\|\hat f\|_{\text{\rm min}}\le (N+1) K c^\alpha=K"$ 
and (2.6) holds.
This proves our claim.

\n Thus we have proved that the ``product map''
$\pi_{1|\Cal X_N} \colon \ \Cal X_N\to \tilde A_1$ is
a complete surjection. Now let us show that this implies
that $\Phi(i,c) \le K_1 c^N$ for all $c\ge 1$. To do that,
consider $u\colon \ \Cal A\to B(H)$ unital homomorphism with
$\|ui\|_{cb}\le c$, let $\tilde u\colon \ \tilde A_c \to
B(H)$ be the canonical extension of $u$. Then we have
$$\|\tilde u\|_{cb(\tilde A_1, B(H))} \le K'' \|\tilde
u\pi_{1|\Cal X_N}\| _{cb(\Cal X_N, B(H))}\tag2.9$$ 
but, for any $j$,  we have, 
if $z=1/c$ 
$$c^j\tilde u\pi_{z|E_j} = \tilde u \pi_{1|E_j}$$ hence
$$\leqalignno{\|\tilde u\pi_{1|E_j}\|_{cb(E_j, B(H))} &=
\|c^j \tilde u\pi_{z|E_j}\|_{cb(E_j, B(H))}&\cr
&\le c^j\|\tilde u\|_{cb(\tilde A_c, B(H))} 
\|\pi_{z|E_j}\|
_{cb(E_j, \tilde A_c)}\cr
\noalign{\text{and by (1.1) and Theorem 1.7}}
&\le c^j.}$$
By Lemma 1.5  this implies
$$\|\tilde
u\pi_{1|\Cal X_N}\| _{cb(\Cal X_N, B(H))}
\le \sum_{j=0}^N \|\tilde u\pi_{1|E_j}\|_{cb(E_j, B(H))}
\le \sum_{j=0}^N c^j\le K'_1 c^N\tag2.10$$
for some constant $K'_1$ independent of $c$. 
Hence (2.9) and (2.10) yield
$$\|\tilde u\|_{cb(\tilde A_1, B(H))}\le K''K'_1c^N.$$
By (2.2) this gives the announced estimate (2.4).

 \n Finally, 
 replacing $g_0,\ldots, g_{N-1}$ by $g'_0,\ldots,
g'_{N-1}$ according to Lemma~1.6 we can
obtain, in the case $1_{\Cal A}\in E$,  an element
$\hat f = \sum^{N-1}_0 z^jg'_j + z^Ng_N$ in $\Cal K\otimes
E^{(N)}$ such that $$\|(I_{\Cal K}\otimes \pi_1)(\hat
f) - f\|_{\text{\rm min}} < 1/2.$$ Moreover by (1.12) and (2.8) we
have 
$$\align
\|\hat f\|_{\text{\rm min}} &\le
\sum\nolimits^{N-1}_0 \|g'_j\|_{\text{\rm min}} + \|g_N\|_{\text{\rm min}}\\
 &\le K_2 = NKc^\alpha\|1_{\Cal A}\|^N + Kc^\alpha.
 \endalign$$
This justifies the last assertion.   
\qed

By a simple modification of the preceding proof we obtain:

\proclaim{Theorem 2.6} Fix a number $\alpha>0$ and let $N$ be an integer
with $N\le \alpha < N+1$. 
Let $X\subset \Cal K$ be a closed subspace for which there is
a projection $P\colon \ \Cal K\to X$ with $\|P\|_{cb}= 1$.
Assume that there is a constant $K$ such that, for any $f$ in
$X\otimes \Cal A$ we have
$$\forall c\ge 1\qquad \|f\|_{X \otimes_{\min} \tilde A_c} \le K c^\alpha
\|f\|_{X \otimes_{\min} \tilde A_1}.\tag2.11$$
 Then the restriction of $I_X\otimes \pi_1$
to $X \otimes_{\min}  OA_N(E)$ is a   surjection 
of $X \otimes_{\min} OA_N(E)$
 onto $X \otimes_{\min} \tilde A_1$, i.e.\
there is a constant $K'$ such that, for any $f$ in
 ${X} \otimes_{\text{\rm
min}} \tilde A_1$ with $\|f\|_{\text{\rm min}} <1$, 
there is $\hat f$ in ${X}
\otimes_{\text{\rm min}} OA_N(E)$ with $\|\hat f\|_{\text{\rm min}} < K'$
 such that
$(I_{X} \otimes \pi_1)(\hat f) = f$.
 Furthermore, this last property
implies that   (2.11) actually holds with $\alpha=N$
 for some (possibly different) constant $K$.   
Finally, if $1_{\Cal A}$ belongs to $E$, then the   
restriction of $I_X\otimes \pi_1$ defines
  a  surjection from $X \otimes_{\min}{E_N}$ 
onto $X \otimes_{\min}\tilde A_1$.
\endproclaim

\proclaim{Corollary 2.7} Consider our usual setting $(E,\Cal A,i)$ and assume 
that
the similarity property holds. Let $d$ be the infimum of the numbers $\alpha>0$ 
for which
there is a constant $K$ such that 
$\Phi(i,c)\le Kc^\alpha$ for all $c$ large enough. Then $d$ is an integer. 
Moreover, there is a
constant $K'$ such that for all 
$c\ge 1$ we have
$$c^d\le \Phi(i,c)\le K' c^d.\tag2.12$$
We will call $d$ the ``similarity degree'' of our setting $(E, \Cal A,i)$.
\endproclaim

\demo{Proof} Let $N$ be the integer such that $N\le d <N+1$.
Then Theorem 2.5 implies $d\le N$, hence $d=N$. Thus $d$ is an integer. Fix 
$\alpha <d$. Then by
Lemma 2.3, we have necessarily
$\Phi(i,c)\ge c^\alpha$ for all $c\ge 1$. By continuity, this must hold also for 
$\alpha=d$,
whence the left side
of (2.12). Finally, the right side
of (2.12) follows from the last part of Theorem 2.5.
\qed
\enddemo

\remark{Remark 2.8} The case $d=0$ is of course trivial,
 this case happens iff $\Cal A$ is one dimensional. 
The
 case $d=1$
also is trivial, although a bit more interesting. 
By Theorem 2.5, $d=1$ happens
only if the operator space $\tilde A_1$ is completely
isomorphic to a quotient space of the direct sum
of $\comp$ with $E$. For instance, in the situation of
the basic 
Example 1.1, we have $d=1$ only if $C^*(G)$ is 
completely isomorphic to a quotient space of $\ell_1 (\Gamma)$,
or equivalently only if $C^*(G)$ is a $\max$-space, in the sense
of [BP1]. By [BP1, Pa5], we know that this can happen only
if $C^*(G)$ is finite dimensional, whence only if
(and a posteriori iff) $G$ is finite.
\endremark

\remark{Remark 2.9} For simplicity, we will identify $E$ with $i(E)$ in
this remark, so we view $E$ simply as a subset of $\Cal A$. We also view
$\Cal A$ as a subset of $\tilde A_1$. We will moreover assume that $E$
contains the unit. Then, by Theorem~2.5, the degree $d$  (as defined in
Corollary~2.7) is equal to the smallest integer $d$ with the property that
the natural product map from $E \otimes_h\cdots \otimes_h E$ ($d$ times)
to $\tilde A_1$ is a complete surjection. By the very definition of the Haagerup 
tensor product, this last
property can be restated as follows:\ there is a constant $K$ such that
for any $n$, any $\vp>0$ and any $a = (a_{ij})$ in $M_n(\Cal A)$ we can find 
matrices
$x_1,\ldots, x_d$ with (say) $x_1\in M_{q_1q_2}(E)$, $x_2\in
M_{q_2q_3}(E),\ldots, x_d\in M_{q_dq_{d+1}}(E)$ and with $q_1 = n =
q_{d+1}$ so that the matricial product
$x_1\cdot x_2\ldots x_d$
 (this is a product in $M_n(\Cal A)$, recall that we view $E$ as
a subset of $\Cal A$) satisfies
$$\|a-x_1\cdot x_2\ldots x_d\|_{M_n(\tilde A_1)}<\vp,$$
 and finally we have
$$\prod^d_{i=1} \|x_i\|_{M_{q_iq_{i+1}}(E)} \le
K\|a||_{M_n(\tilde A_1)}.$$
\endremark

In most of the ``concrete'' examples considered below the
space $E$ is a ``maximal'' operator space in the sense of
[BP1]. In that case, we may apply the decomposition
decribed in (0.4)
to any rectangular matrix $x$ in $M_{pq}(E)$ (by just
adding enough zeros to make it a square matrix). 
Using this
fact, we obtain the following.

\proclaim{Proposition 2.10} Consider a setting
$(i,E,\Cal A)$. Assume that the operator space 
$E$ is a maximal operator
space (in the sense of [BP1]) and
 that $i(E)$ contains
 the unit of $\Cal A$. Then the similarity degree $d$ of $(i,E,\Cal A)$ is
equal to the smallest integer $d$ with the following property:\ there is a
constant $K$ such that for all $n$ any element $x$ in $M_n(\Cal A)$ with
$\|x\|_{M_n(\tilde A_1)} <1$ can be written 
as a limit (in the norm
of ${M_n(\tilde A_1)}$) of  matricial products of the
form (again  we view $E$ as
a subset of $\Cal A$)
$$\alpha_1D_1\alpha_2D_2\ldots D_d\alpha_{d+1}$$
where $\alpha_1,\ldots, \alpha_{d+1}$ are rectangular scalar matrices,
with say $\alpha_i \in M_{p_iq_i}$, $p_1 = n$, $q_{d+1}=n$ and
$D_1,\ldots, D_d$ are diagonal matrices with entries in $E$, with
$D_i\in M_{q_iq_i}(E)$ (with $q_i=p_{i+1}$)
and finally we have
$$\prod^{d+1}_{i=1} \|\alpha_i\|  \prod^d_{i=1} \|D_i\| <
K.$$ (Note that we can assume if we wish, by adding zero
entries, that $q_2 = p_3= q_3 =\cdots= p_d = q_d = N$ for
some $N$ large enough.)
\endproclaim

\head\S 3. Groups\endhead

\n Let $G$ be a discrete group. In this section, we apply our results in the
case
$$E = \Cal A  = \ell_1(G),$$ with $i: E \to A $ equal to the identity.
We equip $E = \ell_1(G)$ with its ``maximal'' operator space structure, so
that for a map $u\colon \ E\to B(H)$ boundedness and complete boundedness
are equivalent and
$$\|u\|_{cb} = \|u\|.$$
Observe that $\Cal A = \ell_1(G)$ is a unital (Banach) algebra for the
convolution product. The unit element of $\Cal A$ is $\delta_e$ defined by
$\delta_e(t) = 1$ if $t=e$ and $0$ otherwise. We have
$\|\delta_e\|_E = 1$.

\n It is classical in this case that
$$\tilde A_1 = C^*(G)$$
the full $C^*$-algebra of $G$. Indeed, any contractive unital homomorphism
$u\colon \ \ell_1(G) \to B(H)$ induces a norm one representation $\pi\colon
\ G\to B(H)$ which is {\it automatically\/} a unitary representation.
(Indeed $\|\pi(g)\|\le 1$ and $\|\pi(g)^{-1}\|\le 1$
implies $\pi(g)$ unitary for any $g$ in $G$.) It also is a
classical fact that the dual of $\tilde A_1 = C^*(G)$ can
be identified with the space $B(G)$ of all coefficients of
the unitary representations of $G$ (\cf [Ey, FTP]). The space $B(G)$ is
defined as the space of all functions $\varphi\colon \
G\to {\comp }$ for which there is a unitary representation
$\pi\colon \ G\to B(H_\pi)$ and vectors $\xi,\eta \in
H_\pi$ such that
 $$\forall\ t\in
G\qquad \varphi(t) = \langle
\pi(t)\xi,\eta\rangle.\tag3.1 $$
 Moreover one defines $$\|\varphi\|_{B(G)} =
\inf\{\|\xi\|\, \|\eta\|\}$$ where the infimum runs over
all possible representations of $\varphi$ as in (3.1). One can
imitate this definition for the algebra $\tilde A_c$:\

\n Let us denote by $B_c(G)$ the space of all functions
$\varphi\colon \ G\to~{\comp }$ for which there is a
uniformly bounded representation $\pi\colon \ G \to
B(H_\pi)$ with $|\pi|\le c$ and vectors $\xi,\eta$ in
$H_\pi$ such that (3.1) holds. We then define again
$$\|\varphi\|_{B_c(G)} = \inf\{\|\xi\|\, \|\eta\|\}\tag3.2$$ 
where the infimum runs over all possible such
decompositions of $\varphi$. If $c=1$, we recover the
unitary case so that $B_1(G)$ is identical to $B(G)$. Now
consider a function $f$ in $\Cal A = \ell_1(G)$. Clearly we
have $$\|f\|_{\tilde A_c} = \sup\left\{\left| \sum_{t\in
G} f(t) \varphi(t)\right|\ \bigg|\ \varphi\in B_c(G)\quad 
\|\varphi\|_{B_c(G)} \le 1 \right\}.\tag3.3$$ 
More
precisely, we have the following well known fact.

\proclaim{Proposition 3.1} With the natural duality appearing in (3.3) we
have
$$B_c(G) = (\tilde A_c)^*$$
with equal norms.
\endproclaim

\demo{Proof} By (3.3) the unit ball of $B_c(G)$ (which is convex) is weak-$*$ 
dense
in the unit ball of $(\tilde A_c)^*$. Hence for any $\varphi$ in the unit
ball of $(\tilde A_c)^*$ there is a net $\varphi_i$ in the unit ball of
$B_c(G)$ which tends pointwise to $\varphi$. Then, by a standard
ultraproduct argument, one can check that $\varphi$ itself is in the unit ball
of $B_c(G)$.\qed
\enddemo

We will also need the space of Herz-Schur multipliers on $G$ which we
denote by $M_0(G)$, we refer to [DCH, BF1-2, Bo1, H3, P2] for more
information. We recall that a function $\varphi\colon \
G\to {\comp }$ is in the space $M_0(G)$ iff there are
bounded Hilbert space valued functions $x\colon \ G\to H$
and $y\colon \ G\to H$ such that $$\varphi(s^{-1}t) =
\langle x(t), y(s)\rangle.\leqno \forall\ s,t\in G$$
Moreover, we denote $$\|\varphi\|_{M_0(G)} =
\inf\{\sup_{s\in G} \|x(s)\|_H \cdot \sup_{t\in
G}\|y(t)\|_H\}\tag3.4$$ 
where the infimum runs over all possible
factorizations of $\varphi$. For the reader's convenience,
we will now reformulate explicitly the meaning of the
constants introduced in the previous section. Let $c\ge 1$.

Consider a bounded representation $\pi\colon \ G\to B(H)$ 
with $|\pi|\le
c$. Assume that $\pi$ is unitarizable then we denote
$$\text{Sim}(\pi) = \inf\{\|S^{-1}\|\, \|S\|\}$$
where the infimum runs over all invertible operators $S\colon \ H\to H$
such that $t\to S^{-1}\pi(t)S$ is a unitary representation. Then we set
$$\Phi_G(c)  = \sup\{\text{Sim}(\pi)\}\tag3.5$$
where the sup runs over all uniformly bounded representations with
$|\pi|\le c$. Let $i_G\colon \ E\to~\Cal A$ be the setting associated to the
identity map of $\ell_1(G) = E=\Cal A$. By Proposition~1.8 and (2.2)', we
know that $$\Phi_G(c) = \Phi(i_G,c)=\|\tilde A_1\to \tilde A_c\|_{cb}.\tag3.5'$$
It will be convenient for our discussion to introduce also
$$\Psi_G(c) = \sup\{\|f\|_{B(G)}\mid f\in B_c(G)
\quad \|f\|_{B_c(G)} \le 1\}.$$ Note that by (2.2)$'$ the
inclusion $\tilde A_1 = C^*(G) \to \tilde A_c$ has norm
$\le \Phi_G(c)$, hence we have by Proposition~3.1
$$\Psi_G(c) \le \Phi_G(c).\tag3.6$$
 Moreover, again
by Proposition 3.1 $$\Psi_G(c) = \|\tilde A_1\to \tilde
A_c\|.$$

\proclaim{Theorem 3.2} The following properties of a discrete group $G$ are
equivalent:\medskip
\roster
\item"{\text{\rm (i)}}" $G$ is amenable.
\item"{\text{\rm (ii)}}" $\Phi_G(c) \le c^2$ and $\Psi_G(c) \le c^2$ for all 
$c>1$.
\item"{\text{\rm (iii)}}" There is $\alpha<3$ and a constant $K$ 
 such that for all $c\ge 1$
$$\Phi_G(c) \le Kc^\alpha.$$
\item"{\text{\rm (iv)}}" $\Phi_G(c) < c^3$ for some $c>1$ and $\Phi_G(c)<\infty$  
 for 
all $c>1$.
\item"{\text{\rm (v)}}" $\exists \alpha <3\  \exists K>0\  $ such that
for all $c\ge 1$
$$\Psi_G(c) \le Kc^\alpha.   $$
\endroster
\endproclaim

\demo{Proof} (i) $\Rightarrow$ (ii) is Dixmier's classical result [Di].
(ii)~$\Rightarrow$~(iii)~$\Rightarrow$~(iv) are trivial by (3.6).
 Moreover,
(iv)~$\Rightarrow$~(iii) follows from Lemmas~2.2 and 2.3
and (3.5)'. In addition,  
(iii)~$\Rightarrow$~(v) follows from (3.6).

\n Hence it remains only to prove (v) $\Rightarrow$ (i).
Assume (v). By Theorem~2.6 with $X=\comp$, the restriction of $\pi_1$ to
$E_2$ is a surjection from $E_2$ onto
$\tilde A_1 = C^*(G)$. Equivalently, this means that the
adjoint map $w_2\colon \ B(G)\to (E_2)^*$ is an
isomorphic embedding, so that for some $\delta>0$ we have
$$\forall\ \varphi\in B(G)\qquad \delta\|\varphi\|_{B(G)} \le 
\|w_2(\varphi)\|_{E^*_2}. \tag3.7
$$ Now
assume $\varphi$ finitely supported. We have
$$\|w_2(\varphi)\|_{E^*_2} =
\sup\left\{\left|\sum_{s,t\in G} \varphi(st)
\alpha(s,t)\right|\right\}$$ where the supremum runs over
all $\alpha = \sum\limits_{s,t\in G} \alpha(s,t)
\delta_s\cdot \delta_t$ in the unit ball of $E_2$.
By Proposition~1.10, the space $E_2$ can
be naturally identified with $\ell_1(G) \otimes_h
\ell_1(G)$, so that for $\alpha$ as above
$$\|\alpha\|_{E_2} = \left\|\sum_{s,t\in G}
\alpha(s,t) \delta_s\otimes
\delta_t\right\|_{\ell_1(G)\otimes _h\ell_1(G)}.$$ Hence
we find $$\|w_2(\varphi)\|_{E^*_2} =
\left\|\sum_{s,t\in G} \varphi(st) e_s\otimes
e_t\right\|_{(\ell_1(G) \otimes_h \ell_1(G))^*}\tag3.8$$
 where $e_t\in \ell_1(G)^*  = \ell_\infty(G)$ is
biorthogonal to $\delta_t$, i.e.\ $e_t(\delta_s) = 1$ if
$t=s$ and 0 otherwise. But now it is well known
 (\cf [DCH, or P1]) that the right side of (3.8) is equal to
$\|\varphi\|_{M_0(G)}$. Hence we deduce from (3.7) that
for all finitely supported $\varphi\colon \ G\to {\comp }$
$$\delta\|\varphi\|_{B(G)} \le \|\varphi\|_{M_0(G)}.$$ By
a result due to Bo\.zejko [Bo2], this implies that $G$
is amenable, whence (i).\qed
\enddemo

\proclaim{Corollary 3.3} If $G$ is not amenable, 
then for any $c>1$ there is a representation
$\pi_c : \ G\to B(H_c)$ with $\|\pi_c\|\le c$ such that 
$c^3\le \inf\{\|S\| \|S^{-1}\|\}$, where the infimum runs
 over all the similarities $S$ such that
$S^{-1}\pi_c(.) S$ is unitarizable. 
\endproclaim

\demo{Proof} By the preceding statement, we know that $\Phi_G(c)\ge c^3$ for all 
$c>1$. 
If
$\Phi_G(c)=\infty$, we clearly have the conclusion. Otherwise 
$\tilde A_1$ and $\tilde A_c$      are isomorphic. Then we represent $\tilde 
A_c$  as a subalgebra
of some
$B(H)$, say
$ \tilde A_c  \subset B(H_c)$ and we define $\pi_c$ to be the
representation on $G$ associated to the restriction to
$E=\ell_1(G)$  of the canonical morphism  $\tilde A_1\to \tilde A_c$. 
By (2.2)' we have  $\|\tilde A_1\to \tilde A_c\|_{cb} \ge c^3 $, hence by
  Proposition 1.8, this representation  has the desired property. \qed
\enddemo

The next result recapitulates what we know from \S 2.

\proclaim{Theorem 3.4} Assume that every uniformly bounded group
representation on $G$ is unitarizable.
\medskip
\roster
\item"{\text{\rm (i)}}" Then the function $\Phi_G(c)$ defined in (3.5) is finite 
for 
all
$c\ge 1$. Moreover, let $d(G)$ be the smallest $\alpha>0$ such that
$\Phi_G(c) \in O(c^\alpha)$ when $\alpha\to \infty$. Then $d(G)$ is an
integer. We call it the similarity degree of $G$.
\item"{\text{\rm (ii)}}" We have $\Phi_G(c) \ge c^{d(G)}$ for all $c\ge 1$ and 
$d(G)$ 
is
the largest integer with this property.
\item"{\text{\rm (iii)}}" The degree $d(G)$ is the smallest integer $N$ such 
that the
natural ``product'' mapping
$$W_N\colon \ \underbrace{\ell_1(G)\otimes _h\cdots \otimes_h
\ell_1(G)}_{N\ \text{\rm times}} \longrightarrow C^*(G)$$
which takes $\delta_{t_1} \otimes \cdots\otimes \delta_{t_N}$ to
$\delta_{t_1t_2\ldots t_N}$ is a complete surjection onto $C^*(G)$.
\endroster
\endproclaim

\demo{Proof} The first part follows from Theorem 2.5 and especially from (2.4). 
By
Lemma~2.3, if $\Phi_G(c) < c^\alpha$ for some $\alpha < d(G)$ then $\Phi(c)
\in O(c^\alpha)$. Therefore we must have $\Phi_G(c) \ge c^\alpha$ for all
$c\ge 1$ and $\alpha < d(G)$, whence the second part.

\n Finally, the third part follows from Theorem~2.5 and
Proposition~1.10, which tell us that $\Phi_G(c) \in O(c^N)$
when $c\to \infty$ iff $W_N$ is a complete surjection.\qed
\enddemo

\remark{Remark 3.5} With the preceding notation, Theorem 3.2 says that $d(G)\le 
2$ iff $G$ is amenable.
\endremark

\remark{Remark 3.6}
We now return to the Example 1.1. Let $G$ be a discrete group. Let $\Cal A$ be 
the group
algebra of $G$, i.e.\ $\Cal A = \ell_1(G)$ equipped with convolution. Let
$\Gamma\subset G$ be a set of generators for $G$ and let
$E=\ell_1(\Gamma)$, equipped again with its natural
 (=maximal) operator space structure.
Here again, we have $\tilde A_1  = C^*(G)$, but the similarity degree
now depends very much on the choice of the generators.
Let us denote by $d=d(\Gamma,G)$ the similarity
 degree for this setting,
according to Corollary 2.7. Then, by Theorem 2.5, the 
product map
$$p: (t_1,...,t_d)\in \Gamma^d \to t_1...t_d\in G$$
extends to a (complete) surjection from
$\ell_1(\Gamma)\otimes_h...\otimes_h \ell_1(\Gamma)$
onto $C^*(G)$. Let us denote by $[\Gamma]^d\subset G$
 the image
of $\Gamma^d$ under $p$. Then,  the elements
supported by $[\Gamma]^d$ must be dense in $C^*(G)$,
and a fortiori, say, in $\ell_2(G)$. This clearly
implies  (denoting by $e$ the unit element of $G$)
$\{e\}\cup \cup_{j\le d}[\Gamma]^j=G.$
Therefore, every element of $G$ can be written
 as a product of at most
$d$ elements of $\Gamma$. If we introduce
 the usual distance
on $G$ relative to (the Cayley graph of)
 $\Gamma$, this means that the diameter of $G$ is at most $d$.
This remark allows to produce examples of similarity settings
with arbitrarily large finite
similarity degrees. Indeed, just consider 
$G=\ent^N$, for some integer $N$, and take for $\Gamma$ the subset formed of all 
elements
with only one non zero coordinate. Clearly, by the preceding
remarks we have $N\le d$ in this case. In the converse
 direction,
we claim that $d\le 2N$. Indeed, let
$\pi : G\to B(H)$ be a representation such that
$\sup_{t\in \Gamma}\| \pi(t)\|\le c$. Then clearly 
   $\sup_{t\in G}\| \pi(t)\|\le c^N$. Now, since $G$ is amenable,
this implies by Dixmier's theorem that $\text{\rm Sim}(\pi)\le c^{2N}$.
Hence, we have shown that $N\le d\le 2N$. 
\endremark

More precisely, consider in $C^*({\ent}^N) \simeq C^*({\ent})
\otimes_{\text{\rm  min}} \cdots \otimes_{\text{\rm min}} C^*({\ent})$ the
subspace $E = C_1 +\cdots+ C_N$ with $C_i  = 1\otimes\cdots\otimes1 \otimes
C^*({\ent}) \otimes 1\otimes\cdots\otimes1$
 where $C^*({\ent})$ appears at the $i$-th place.
We equip $E$ again with its maximal operator 
space structure and we let $\Cal A$ 
be the algebra generated by $E$ in $C^*({\ent}^N)$. We will show that
 the  degree of the ``similarity setting" constituted of 
the inclusion $E\subset \Cal A\subset C^*({\ent}^N)$ is equal to $2N$.

\n Let $\pi\colon \ C^*({\ent}^N)\to B(H)$ be a unital homomorphism
such that
$\|\pi_{|E}\|\le c$. Then clearly $\sup_{t\in {\ent}^N}
\|\pi(t)\|\le c^N$, hence by Dixmier's Theorem $\|\pi\|_{cb} \le c^{2N}$.
On the other hand, since $d(C^*({\ent}))=2$, there exists a unital
homomorphism $u \colon\  C^*({\ent})\to B(H)$ with $\|u\|\le c$ and
$\|u\|_{cb}\ge c^2$. Let $\pi = u\otimes u \otimes\cdots\otimes u\colon \
C^*({\ent}^N)\to$ $B(H)\otimes_{\min}\cdots\otimes_{\min} B(H)$. Then $\pi$ is a
unital homomorphism and it is not hard to check that $\|\pi_{|E}\|\le 1+2Nc$.

\n Indeed, if $x = \sum\limits^N_1 x_i\in E$ with $\|x\|\le 1$ and if
$\tau(x_i)$ is the usual trace, i.e.\ the Haar integral on 
$C(\T^N)=C^*({\ent}^N)$, then
$ \left|\sum^N_1 \tau(x_i)\right|\le 1$ and $\left\|\sum^N_1
(x_i-\tau(x_i))\right\| \le 2$, which implies (using conditional expectations)
$\|x_i-\tau(x_i)\| \le 2$, hence
$$\align
\|\pi(x)\| &= \left\|\sum\nolimits^N_1 1\otimes
u(x_i)\otimes 1\right\|\\
&\le \left\|\left(\sum\nolimits^N_1 \tau(x_i)\right)
1\otimes\cdots\otimes1\right\| + \left\|\sum\nolimits^N_1 1\otimes\cdots
(u(x_i)-\tau(x_i)1)\otimes\cdots\otimes1\right\|\\
&\le 1 + \sum\nolimits^N_1 \left\|u(x_i-\tau(x_i)1)\right\| \le 1 +
Nc\|x_i-\tau(x_i)1\|\\
&\le 1 + 2Nc.
\endalign$$
On the other hand, we clearly have $\|\pi\|_{cb} \ge (c^2)^N = c^{2N}$,
hence this proves that the degree $d$ of this similarity setting is {\it
exactly\/} equal to $2N$.

\remark{Remark} In the   setting described in Example 1.1, 
let $G$ be a discrete amenable group, 
so that $\tilde A_1 =
C^* (G)= C^*_\lambda(G)$. We claim that the smallest constant $K$ appearing in
Proposition 2.10 (with $d=2$) is actually equal to 1.
\endremark

\n Indeed, consider $x$ in $M_n(C^*_\lambda(G))$, with
$\|x\|<1$. Since $G$ is amenable, its Fourier algebra
$A(G)$  has an approximate unit in its unit ball. Hence, we
may assume (by density) that $x$ is of the following form
$x = \sum\limits_{t\in G} y(t) \otimes \lambda(t)
\varphi(t)$ with $y(t) \in M_n$ and $y = \sum\limits_{t\in
G} y(t) \otimes \lambda(t)$ such that
$\|y\|_{M_n(C^*_\lambda(G))} <1$ and $\varphi$ of the form
$\varphi(t) = \langle\lambda(t) \xi,\eta\rangle$ with
$\|\xi\|\, \|\eta\|<1$ where $\xi(\cdot), \eta(\cdot)$,
$y(\cdot)$ and $\varphi(\cdot)$ are all finitely supported.

\n Equivalently we have
$$x_{ij} = \sum_{s,\theta} y_{ij}(s\theta^{-1}) \lambda(s)
\lambda(\theta^{-1}) \overline{\eta(s)} \xi(\theta).\tag3.9$$
Then a simple computation shows that we can write
$$x = A_1D_1A_2D_2A_3\tag3.10$$
where $A_1,A_2,A_3$ are rectangular scalar matrices and where $D_1,D_2$
are diagonal with entries in $A$ of the form $\lambda(t)$ for some $t$ in
$G$. Moreover, we have 
$$\|A_1\|\, \|A_2\|\, \|A_3\|<1. \tag3.11$$
Explicitly, we can take
$$\align
&A_1(i,(s,k)) = \overline{\eta(s)} \delta_{ik}\\
&A_2((s,k), (\theta,\ell)) = y_{k\ell}(s\theta^{-1})\\
&A_3((\theta,\ell),j) = \xi(\theta)\delta_{\ell j}
\endalign$$
and the diagonal matrices defined by
$$\align
&D_1((s,k), (s,k)) = \lambda(s)\\
&D_2((\theta,\ell), (\theta,\ell)) = \lambda(\theta^{-1}).
\endalign$$
Note that we can restrict the sums in (3.9) to be over the finite
subsets of $G$ where $\xi$ and $\eta$ are supported, so that we indeed
obtain finite matrices in (3.10),  and (3.11) is easy to
check.

\n Thus the decomposition (3.10) clearly implies our claim
that Proposition 2.10 holds with $d=2$ and $K=1$.

\head \S 4. Operator algebras\endhead

\n We now come to the main application of our results. Let $A$ be a unital 
operator algebra. Let
$$E = \max(A)$$
in the sense of [BP1] (see (0.4) above). The operator space $E$ is equal to $A$ 
as a Banach 
space, but its operator space structure is characterized by the property 
that, for any linear map $u\colon \ E\to B(H)$, we have
$$\|u\| = \|u\|_{cb}.\tag4.1$$
Here, of course we take $\Cal A=A$,
 and we   let $i\colon \ E\to A$ be the identity of $A$. Of course,
we have $\tilde A_1 =A$ isometrically
(but perhaps not completely so).

\n Then we denote by $d(A)$ the similarity degree of this setting $(i,E,A)$. 

\n Note that, by definition,  $d(A)\le \alpha$ iff there
is a constant $K$ such that, for any bounded
 unital homomorphism $u\colon\ A\to B(H)$, 
there is an invertible $S$ for which 
$a\to S^{-1} u(a) S$ is contractive and such that
$\|S\| \|S^{-1}\|\le K \|u\|^\alpha.$

It is easy to check that for any closed two sided ideal
$I\subset A$, the quotient space $A/I$ (which, by [BRS],
is an operator algebra) satisfies
 $$d(A/I)\le d(A).$$
Moreover, if $B$ is another unital operator algebra
and if $A\oplus B$ denotes the direct sum
(equipped with the norm $\|(x,y)\|=\max\{\|x\|,\|y\|\}$
and the obvious ``block diagonal" operator algebra structure)
then we have
$$d(A\oplus B)\le \max\{d(A),d(B)\}.$$

Now assume that every unital contractive homomorphism $u\colon 
\ A\to B(H)$ is completely bounded. Then clearly
 $A\simeq \tilde A_1$
completely  isomorphically, and there is a constant $K$ such that
$\|u\|_{cb}\le K$ for all unital contractive homomorphisms $u\colon 
\ A\to B(H)$. This implies that,  
if we define
$$\Phi_A(c) = \sup \{\|u\|_{cb}\} $$
where the supremum runs over all unital homomorphisms $u\colon 
\ A\to B(H)$ with $\|u\|\le c$, then in the present setting we have
$$\forall c\ge 1\quad \Phi(i,c)\le \Phi_A(c)\le K
\Phi(i,c).\tag4.2$$ 
Thus, to recapitulate,
 we  obtain the following two statements
(note that the equivalence between
 (a) and (b) below is due to Paulsen [Pa4]).

\proclaim{Theorem 4.1} Let $A$ be any unital 
operator algebra, then the 
following are equivalent: 
\roster
\item"{\text{\rm (a)}}" Every bounded unital homomorphism $u\colon \
A\to B(H)$ is  similar to a completely contractive one.
\item"{\text{\rm (b)}}" Every bounded unital homomorphism 
$u\colon \ A\to B(H)$ is 
completely bounded.
\item"{\text{\rm (c)}}" There is an integer $d\ge0$ such that, for
some constant $K$,  any unital
homomorphism $u\colon \ A\to B(H)$ satisfies $\|u\|_{cb}
\le K\|u\|^d$.
\endroster
\endproclaim

\demo{Proof}  For
the equivalence between (a) and (b)   (due to Paulsen
[Pa4]), see the above Proposition 1.8. If (a) or (b)
holds, then in  the present setting we have $A\simeq
\tilde A_1$ completely  isomorphically and (4.2) holds.
Thus (b)~$\Rightarrow$~(c) follows from Corollary~2.4, 
and the converse is obvious.\qed
\enddemo

\proclaim{Theorem 4.2}
  For any fixed integer $d\ge0$, the
following properties of a unital 
operator algebra $A$ are equivalent:
\roster
  \item"{\text{\rm (i)}}" There is a  constant $K$ such that  any
unital homomorphism $u\colon \ A\to B(H)$ satisfies
$\|u\|_{cb} \le K\|u\|^d$.
\item"{\text{\rm (ii)}}" There is a number $\alpha $ with $d\le
\alpha<d+1$ for which there exists a  constant $K$ such
that  any unital homomorphism $u\colon \ A\to B(H)$
satisfies   $\|u\|_{cb} \le K\|u\|^\alpha$.
\item"{\text{\rm (iii)}}" The natural
product  mapping $T_d$ from $\max(A)
\otimes_h\cdots\otimes_h \max(A)$ ($d$ times)  onto $A$ is
a complete quotient map, i.e.\ it induces a complete 
isomorphism from the
quotient space $\max(A)\otimes_h\cdots\otimes_h 
\max(A)/\text{ker}(T_d)$ onto $A$.  
\item"{\text{\rm (iv)}}" There
is a constant $K$ such that   any  bounded linear map
$u\colon \ A\to B(H)$ satisfies $\|u\|_{cb} \le
K\|uT_d\|_{cb}$.  
\item"{\text{\rm (v)}}" There is  a
constant $K$ such that  the  following holds:\ assume that
a linear map $u\colon \ A\to B(H)$ is such  that there are
bounded linear maps $v_i\colon \ A\to B(H)$ such that 
$\forall x_1,\ldots, x_d\in A$ $$u(x_1x_2\ldots x_d) =
v_1(x_1) v_2(x_2) \ldots v_d(x_d),$$ then we have
$$\|u\|_{cb} \le K\prod^d_{i=1} \|v_i\|.$$ 
  \item"{\text{\rm (vi)}}" There is  a
constant $K$ such that  the  following holds:
 for all $n$, any element $x$ in $M_n({A})$ with
$\|x\|_{M_n(A)} <1$ can be written, for some
integer $N$, as a matricial product of
the form
$$x = \alpha_0D_1\alpha_2D_2\ldots D_d\alpha_{d}$$
where $\alpha_0\in M_{n N}$, $\alpha_1\in M_{N}$,..., $\alpha_{d-1}\in M_{N}$,
$\alpha_{d}\in M_{Nn}$ are scalar matrices (\ie $\alpha_0$ and $\alpha_{d}$
are  rectangular of size $n\times N$ and $N\times n$, and the others are square
matrices of size $N\times N$),
 and
$D_1,\ldots, D_d$ are  $N\times N$ diagonal matrices with entries in $A$,
and finally we have
$$\prod^{d}_{i=0} \|\alpha_i\|  \prod^d_{i=1} \|D_i\| <
K.$$
\item"{\text{\rm (vii)}}" There is a constant
$K$ such that any $x$ in the unit ball of $\Cal K\otimes_{\text{\rm min}} A$
can be written as a product of the form
(recall that
$C_0 \subset \Cal K$ denotes the subspace of diagonal
operators)
$$x = \alpha_0 D_1\alpha_1D_2\ldots
D_d\alpha_{d}$$ with $\alpha_i \in \Cal
K\otimes 1$ and $D_i \in C_0 \otimes_{\text{\rm min}}
A$ such that
$$\prod^{d}_{i=0} \|\alpha_i\| \prod^d_{i=1} \|D_i\| \le K.$$
\endroster
\endproclaim

\demo{Proof}    (i)~$\Rightarrow$~(ii) is trivial.
(ii)~$\Rightarrow$~(iii) follows from Theorem~2.5 and (iv)
is merely a  restatement of (iii). Similarly
(iv)~$\Rightarrow$~(v) is clear by the  known
factorization property of c.b.\ multilinear maps ([CS1-2,
PaS]) and by  (4.1). 

\n  Furthermore, if (v) holds, then every
bounded unital homomorphism  $u\colon \ A\to B(H)$ satifies
$$\|u\|_{cb} \le K \|u\|^d,$$
whence (i). This proves the equivalence of (i)-(v).
 The equivalence
between (vi) and (vii) can be checked by a standard
argument left to the reader.

\n Finally, note that (vi)~$\Rightarrow$~(v) is obvious (with the same constant
$K$). It remains   to prove the converse, so assume (v), then Proposition 2.10
is of course valid in the present setting 
($E=\max(A)$, and $\Cal A=A$) with degree $d$. 
Note that $\|x\|_A= \|x\|_{\tilde A_1}$ for any $x$ in $A$,
 hence for any $a$ in $M_n(A)=M_n(\tilde A_1)$, we have obviously
$$\| a\|_{M_n(\max(A))}\le n^2 \| a\|_{M_n(\tilde A_1)}.$$ This implies, by
(0.4),
 that  any $a$ with $\| a\|_{M_n(\tilde
A_1)}< n^{-2}$can be written as $a=\alpha D \beta$ as in (0.4).
Thus, by Proposition 2.10, any $x$ with  $\| x\|_{M_n(A)}< 1$
can be written as a sum $x=a+ y$ with $y=\alpha_1D_1\alpha_2D_2\ldots
D_d\alpha_{d+1}$ having factors  $\alpha_1,D_1,\alpha_2,D_2,\ldots
,D_d,\alpha_{d+1}$   as in (vi), and with $a=\alpha D \beta$ as above. Now,
by adding redundant factors equal to the unit, we can assume that $a$ is of the
same form as $y$, say  $a=\alpha'_1D'_1\alpha'_2D_2\ldots D'_d\alpha'_{d+1}$,
and then changing $N$ to $2N$ (and $K$ to $K+1$), it is easy to rewrite the sum
$x=a+ y$ as a single product as in (vi). This shows that (v) implies (vi) and
concludes the proof.
 \qed
\enddemo
 
\proclaim{Corollary 4.3} Assume (4.2). Then the similarity degree $d(A)$ 
(defined
as $d(A)=d(i,E,A)$ with $E=\max(A)$ and $i=I_A$) is equal 
to the smallest integer $d$
satisfying the 
equivalent conditions \text{\rm (i)-(vi)} in Theorem 4.2.
 Moreover, when this integer is finite, for any $c>1$, there is
a unital homomorphism $u_c\colon \ A\to B(H)$ with 
$\|u_c\|\le c$ such that
$$\|u_c\|_{cb}\ge c^{d(A)}.$$
\endproclaim

\demo{Proof} The first assertion is clear. The second one follows from
(2.2)', (2.12)  and the obvious fact that, for any unital operator
algebra $A$, in the present setting the natural inclusion
of $\tilde A_1$ into $A$ is completely contractive. \qed 
\enddemo

\remark{Remark 4.4} As in Remark 2.8, $d(A)=1$ iff
$\tilde A_1=\max(A)$   completely isomorphically.
 In contrast with the
group case (see Remark 2.8), there are
infinite dimensional examples when this happens. Indeed,
consider a closed infinite dimensional  subspace $E\subset
B(H)$, which is a maximal operator space, \ie such that
$E=\max(E)$.  Then consider (as in [Pa5]) the subalgebra $A_E\subset
B(H\oplus H)$ formed of all elements of the form
$\left(\matrix
\lambda I&x\\ 0&\lambda I\endmatrix\right)$ with
$x\in E$ and $\lambda\in \comp$. Clearly,
$A_E\simeq \comp\oplus E$, hence $A_E\simeq\max(A_E)$
completely isomorphically, and (4.2) holds in this case, so this gives an 
example with
degree 1, \ie we have $d(A_E)=1$. For examples with
degree equal to 2 and 3, see \S 6 below. We conjecture
that  the value of $d(A)$ can be any integer, but, at the
time of this writing, we do not have any example with $3<
d(A) <\infty$. 
\endremark

\remark{Remark 4.5} In the present setting $(i,E,A)$ as defined in the
beginning of \S 4, the similarity property holds iff we have:
\roster
\item"{(SP)}" Every bounded unital homomorphism $u\colon \ A\to B(H)$ is
similar to a contractive one.
\endroster
\endremark

\n Clearly, by \S 2, this holds iff in this setting the
degree is finite. Up to now, in this section, we have concentrated on algebras 
$A$ which satisfy (4.2). Nevertheless, the above property (SP) could
be of interest even if the right side of (4.2) fails. Note
however that, if we replace $A$ by $\tilde A_1$, then we
return to the situation discussed in Theorem~4.2. More
precisely, the setting being still the same as throughout
this section, we have  $$\Phi_{\tilde A_1}(c) =
\Phi(i,c),\leqno \forall c\ge 1$$ 
and $A$ satisfies (SP)
iff every bounded unital homomorphism $u\colon \ \tilde
A_1 \to B(H)$ is c.b. Note that here $A$ and $\tilde A_1$
are isometric, but perhaps not completely isomorphic.

\remark{Remark 4.6} Fix $n\ge 1$. Let $\Cal U$ be the
unitary group
in $M_n$ with normalized Haar measure $m$. It is not hard
to show that the mapping
$$x\to \int_{\Cal U} xu\otimes u^* \ m(du)$$
is completely contractive   from $M_n$ to 
$\max(M_n)\otimes_h \max(M_n)$.
\endremark

\n Thus, in the case $A=M_n$, the surjection appearing in
Theorem 4.2 (iii) (with $d=2$ here) actually admits a
completely contractive lifting. Consequently,
 when $A=\Cal K$
and $d=2$, the constant $K$ appearing 
in (iii),(iv) (v) or (vi) 
in Theorem 4.2 is actually equal to $1$.
Probably a more general
result holds in the context of
``normal virtual diagonals" in the sense of [E].

\remark{Remark} It is probably possible to develop
 the theory
of the ``similarity degree" in the category of {\it dual}
operator algebras, replacing the Haagerup
tensor product by the   dual variant considered
in [BS] and restricting attention
to weak-$*$ continuous homomorphisms, 
 but we have not pursued this yet. (Note added
may 97: this program has now been successfully
carried out by C. Le Merdy.)
\endremark

\remark{Remark 4.7} Recently, Kirchberg [Ki] showed that a
unital $C^*$-algebra $A$ has the similarity property (in
other words $d(A)<\infty$) iff every derivation
$\delta\colon \ A\to B(H)$, relative to an arbitrary
$*$-representation $\pi\colon \ A\to B(H)$ (we will call
such a derivation   a $\pi$-derivation) is inner. 
Equivalently, we have $d(A)<\infty$ iff there is a
constant $K$ such that any such derivation $\delta$
satisfies $$\|\delta\|_{cb} \le K\|\delta\|.\tag4.3$$
More precisely, a simple adaptation of Kirchberg's
argument shows that (4.3) implies $$d(A)
\le K.\tag4.4$$
 Here is a brief sketch:\ we follow the
presentation of Kirchberg's argument in [P1, p.~129]. Let
$\pi\colon \ A\to B(H)$ be a unital $*$-representation and
let $S\colon \ H\to H$ be self-adjoint and invertible. Let
$\pi_S(x) = S^{-1}\pi(x)S$ and let $\delta(x) =
\text{Log}(S) \pi(x) - \pi(x) \text{Log}(S)$ ($x\in A$).
We assume $\|\pi_S\|=c$ and $\|\pi_S\|_{cb} = \|S^{-1}\|\,
\|S\|$. Fix $a$ unitary in $A$. Consider the entire
function $f(z) = S^z \pi(a)S^{-z}$. We have $\|f(z)\|\le
c$ if $\text{Re}(z) =
 1$  and $\|f(z)\| = 1$ if
$\text{Re}(z)=0$. Hence by log-subharmonicity, we have
$\|f(z)\| \le c^{\theta}$ if $\text{Re}(z) = \theta$,
$0<\theta<1$. Since $ f(0)=\pi(a)   $ is unitary,
we have $\|f(\theta)\|= \| 1+ \theta
f'(0)\pi(a)^{-1}+o(\theta)\| \le 1+  \theta\text{Log}(c)
+o(\theta)$ when $\theta>0$ tends to zero.
Therefore the Hermitian operator $T=f'(0)\pi(a)^{-1}=
\text{Log}(S)   - \pi(a) \text{Log}(S) \pi(a)^{-1}$
satisfies, for any $h$ in the unit sphere of $H$
$$\langle Th,h \rangle \le  \text{ Log } c.$$
Applying this last estimate with $a$ replaced by its
inverse, we obtain  for any $h$ in the unit sphere of $H$
$$-\langle Th,h \rangle \le  \text{ Log } c.$$
Consequently  $ \|T\|\le
\text{Log}(c).$ Hence, we find $\|\delta(a)\| =
\|f'(0)\|=\|T\|\le  \text{ Log } c$, so that $\|\deltaÊ\|Ê\le
   \text{ Log } c$, whence by (4.3)
$\|\delta\|_{cb} \le K\|\delta\| \le K \text{ Log } c$. By
following Kirchberg's argument as presented in [P1,
p.~130], we then conclude that $\text{Log}\|\pi_S\|_{cb} \le
\|\delta\|_{cb}\le K \text{ Log } c$, hence finally
$\|\pi_S\|_{cb} \le \|\pi_S\|^K$.\qed
\endremark

Let $A$ be a unital operator algebra and let $K(A)$ be the smallest
constant $K$ such that $\|\delta\|_{cb}\le K\|\delta\|$
for any  completely contractive unital homomorphism
$\pi\colon \ A\to B(H)$ and any $\pi$-derivation
$\delta\colon \ A\to B(H)$.

\n Curiously, in almost all of E. Christensen's works in the $C^*$-case,
  the upper estimates   which appear  for $K(A)$,  are
all natural integers (\cf  [C1-4]).  On the other hand,
note that if  $A$ satisfies  (iv)-(vi) in Theorem 4.2, then
we have $K(A)\le K d$. This suggests various questions
which we could not answer (we ask this
   for $C^*$-algebras only, but the questions  make sense in
general):

 \proclaim{Problem 4.8} Is any of the best constants
 $K$ appearing in 
the conditions  (i), or (iv)-(vi)   (from Theorem 4.2)
automatically equal to $1$? 
\endproclaim

\proclaim{Problem 4.9} Is $K(A)$ always an integer when it is finite?
\endproclaim

\head \S 5. Uniform algebras\endhead

\n A uniform algebra is a closed unital subalgebra $A$ of a commutative unital
$C^*$-algebra $C$, such that $A$ generates $C$ as a
$C^*$-algebra.
Equivalently, we can view $A$ as a unital subalgebra of the algebra $C(T)$
of all continuous complex functions on some compact set $T$, which separates the
points of $T$. We say that $A$ is proper if $A\ne C(T)$. A typical example
is the disc algebra $A(D)$ formed of all continuous complex valued
functions on $\partial D$ which extend continuously and analytically inside
$D$. Equivalently, $A(D)\subset C(\partial D)$ can be  viewed as the
closure in $C(\partial D)$ of the space of all polynomials.

Recently, we produced the first example of a bounded unital homomorphism on
$A(D)$ which is not c.b.\ (\cf \ [P5]). 
It is possible that every proper
uniform algebra admits such homomorphisms 
and has infinite degree
(note that the extension to other domains 
of $\comp^n$ such as the polydisc or
the ball is trivial). However, at the time of this writing, the
only general result we have in this direction is the following one.

\proclaim{Theorem 5.1} Let $A$ be a  uniform algebra, such that
any contractive unital homomorphism from $A$ to $B(H)$ is $c.b.$ 
Then $d(A) \ge 3$ if $A$ is proper.
 In other words, $d(A)\le 2$ iff $A$ is a commutative
$C^*$-algebra.
\endproclaim

\remark{Remark 5.2} Equivalently (by [Sh]), $d(A)\le 2$ iff $A$ is an amenable 
Banach
algebra in the sense of e.g.\ [Pi2].
 Compare with Remark 3.5.
 Note that there seem to be no known example 
of an amenable
operator algebra which is not a $C^*$-algebra (see [CL]).
\endremark

The proof is based on the following two results.
To state the first one, it will be convenient to work with a slightly 
unconventional
version of the space $H^\infty$, which we now introduce.
Let ${\bold T}=\partial D$. Consider   $\Omega =
{\bold T}^I$ with
$I =
\{1,2,3,\ldots\}$. Let $(z_i)_{i\ge 1}$ denote the coordinates on $\Omega$
and let $\Cal A_n$ be the $\sigma$-algebra generated by $(z_1,\ldots,
z_n)$ with $\Cal A_0$ the trivial $\sigma$-algebra.
	Let $m$ be the usual
probability measure on ${\bold T}^I$ (= normalized Haar measure). 
 Every $m$-integrable 
function $f\colon \
\Omega\to {\comp }$ defines a martingale $(f_n)_n$ by setting $f_n = {\ee 
}(f\mid \Cal A_n)$.
A martingale $(f_n)_n$, relative to the filtration 
$(\Cal A_n)$, is called ``Hardy'' if for
each
$n\ge 1$ the function $f_n$ depends analytically on $z_n$ (but arbitrarily on
$z_1,\ldots, z_{n-1}$).
 We denote by   $H^{\infty}_{\text{\rm m}}$ the
subspace of
  $L^{\infty}(\Omega,m)$ formed by all $f$ which generate a Hardy
martingale.  

\n In Harmonic Analysis terms, the space  $H^{\infty}_{\text{\rm m}}$  is
indeed the version of  $H^{\infty}$  associated to the
ordered group ${\ent}^{(I)}$ (formed of all the finitely supported families
$n=(n_i)_{i\in I}$ with
$n_i\in {\ent}$), ordered by the {\it lexicographic\/} order, \ie the
order defined by setting $n'<n''$ iff  the last differing coordinate (=``letter" 
with reversed
alphabetical order) satisfies $n'_i<n''_i$.
 As explained e.g.\ in [Ru,
Chapter 8], this group has a ``{\it linear\/}'' behaviour and the associated
$H^p$ spaces on it behave like the classical (unidimensional) ones.
More generally, for any Banach space $X$, we will
denote by $H^p_{\text{\rm m}}(X)$ $(1\le p\le \infty)$ the usual $H^p$-space of 
$X$-valued
functions on the  group $\Omega$ (with ordered dual ${\ent}^{(I)}$), in 
Bochner's sense. 

\proclaim{Lemma 5.3} Let $I$ be any set, and let $X = (\ell_1(I) \otimes_h
\ell_1(I))^*$. Then there is a constant $C$ such that for any Hardy
martingale $(f_n)$ in $H^\infty_m(X)$ we have
$$\sup_N \left\|\sum\nolimits^N_0 \int \bar z_n df_n\, dm \otimes
e_n\right\|_{X\otimes_{\text{\rm min}} \max(\ell_2)} \le 
C\|f\|_{H^\infty_m(X)}.\tag5.1$$
\endproclaim

\demo{Proof} We follow [P2, \S 4]. First observe that it
suffices to prove this when $I$ is a finite set. Indeed, if we know (5.1)
for all finite sets then we can obtain it for an arbitrary set $\Cal I$
by taking the supremum of each side over all finite subsets $I\subset \Cal
I$. Let us assume that $I$ is finite, so that $X = \ell_\infty(I)
\otimes_h \ell_\infty(I)$.

\n It clearly suffices to prove that for all functions $f$ in the open unit
ball of $H^\infty_m(X)$ we have
$$\left\|\sum\nolimits^N_0 \int \bar z_n df_n dm\otimes
e_n\right\|_{X\otimes_{\text{\rm min~max}}(\ell_2)} \le C$$
for some absolute constant $C$ independent of $N$ and $f$. Let $f\in
H^\infty_m(X)$ with $\|f\|_{H^\infty_m(X)}<~1$. 
By invoking [P2, Theorem~4.2
and Remark~4.4], it follows that we can write for $n=0,\ldots, N$
$$\int\bar z_n df_n dm \left(= \int \bar z_n f\, dm\right) = A_n+B_n\tag5.2$$
with $A_n,B_n \in\ell_\infty(I) \otimes_h \ell_\infty(I)$ such that
$$
 \left\|\sum\nolimits^N_0 A_n\otimes e_n\right\|_{\ell_\infty (I,\ell_2)
\widehat\otimes \ell_\infty(I)} + \left\|\sum\nolimits^N_0 B_n\otimes
e_n\right\|_{\ell_\infty(I)\widehat\otimes \ell_\infty(I,\ell_2)} \le
C \tag5.3$$
for some absolute constant $C$. Here we denoted by $\ell_\infty(I,\ell_2)$
the Banach space of all bounded $\ell_2$-valued functions on $I$ equipped
with its natural norm. In addition, we denoted by $\widehat\otimes$ the
projective tensor product and we made the obvious identifications of
$\ell_\infty(I) \otimes \ell_\infty (I)\otimes \ell_2$ with a subset
respectively of $\ell_\infty(I,\ell_2) \otimes\ell_\infty(I)$ and
$\ell_\infty(I) \otimes \ell_\infty(I,\ell_2)$. By a simple argument, one
can check that
$$\left\|\sum\nolimits^N_0 A_n\otimes e_n\right\|_{X\otimes_{\text{\rm 
min~max}}(\ell_2)} \le \left\|\sum\limits\nolimits^N_0 A_n\otimes
e_n\right\|_{\ell_\infty (I,\ell_2) \widehat\otimes \ell_\infty(I)}$$
and similarly
$$\left\|\sum\nolimits^N_0 B_n\otimes e_n\right\|_{X \otimes_{\text{\rm 
min~max}}(\ell_2)} \le \left\|\sum\limits\nolimits^N_0 B_n\otimes
e_n\right\|_{\ell_\infty(I) \widehat\otimes \ell_\infty (I,\ell_2)}.$$ 
  Therefore, we conclude from (5.2) and (5.3) that
$$\left\|\sum\nolimits^N_0 \int \bar z_n df_n dm \otimes
e_n\right\|_{X\otimes _{\text{\rm min~max}}(\ell_2)} \le C.$$\qed
\enddemo

\remark{Remark} It is possible to complete the proof without appealing to
the projective tensor product, remaining in the framework of the Haagerup
tensor product, but this option would unnecessarily lengthen the argument.
\endremark

In [Kis], S.~Kisliakov proved the remarkable fact that, if $A\subset C(T)$ is
any proper uniform algebra, there is no bounded linear projection from
$C(T)$ onto $A$. In [Ga], Garling extended Kisliakov's result. In
particular, the following result is implicit in Garling's paper, but is
proved there (see the proof of Theorem~2 in [Ga]).

\proclaim{Lemma 5.4} Let $A$ be any proper uniform algebra. Then for some
$\beta>0$ there is, for each integer $n$, a Hardy martingale $f_1,\ldots,
f_n$ with values in the unit ball of $A^*$ but such that
$$\left\|\int \bar z_k df_k \right\|_{A^*} \ge \beta \quad
\text{for}\quad k=1,2,\ldots, n.$$
\endproclaim\

\demo{Proof of Theorem 5.1} Let $A\subset C(T)$ be a subspace with the
induced operator space structure. By a joint result due
independently to Junge and to Paulsen and the
author (\cf  [Pa6, Theorem 4.1]) there is a constant $C'$ such that for any
sequence $(\xi_n)$ in $A^*$ we have
$$\left(\sum \|\xi_n\|^2\right)^{1/2} \le C'  \sup_N\left\|
\sum\nolimits^N_0 \xi_n \otimes e_n\right\|_{A^*\otimes_{\text{\rm min}} \max
(\ell_2)}. \tag5.4$$
By Theorem 4.1, if $d(A)\le 2$, then $A$ is a quotient (as an operator
space) of $\max(A) \otimes_h \max(A)$, or a fortiori of $\ell_1(I) \otimes
_h \ell_1(I)$ for some index set $I$. Therefore, there is a subspace
$Y\subset X$ and an complete isomorphism $w\colon\ A^*\to Y$ with
$\|w\|_{cb}\le 1$. By Lemma~5.3 and by (5.4), this implies that for all
Hardy martingales $f=(f_n)$ in $H^\infty_m(A^*)$ we have
$$\left(\sum \left\|\int \bar z_n df_n dm\right\|^2_{A^*}\right)^{1/2} \le
C''\|f\|_{H^\infty_m(A^*)}\tag5.5$$
with $C'' = CC'\|w^{-1}\|_{cb}$.
Finally, by Lemma 5.4, $A$ cannot be proper (since (5.5) would imply
$\beta\sqrt n\le C''$ for all $n$).
\enddemo

\remark{Remark 5.5} The preceding proof establishes more than claimed in
Theorem~5.1. Indeed, we conclude that if $A$ is proper then the operator
space $A^*$ is not completely isomorphic to any subspace of
$(\ell_1(I)\otimes_h \ell_1(I))^*$ for any set $I$. Stated in that form the
result cannot be improved much. Indeed, it can be shown that if $A$ is an
arbitrary operator space, then for a suitable set $I$, $A^*$ embeds
completely isometrically into $(\ell_1(I) \otimes_h \ell_1(I) \otimes_h
\ell_1(I) \otimes_h \ell_1(I))^*$.
\n Indeed, let
$$\Cal X_N(I) = \ell_1(I) \otimes_h\cdots \otimes_h \ell_1(I)\qquad (N 
\text{ times}).$$
Note that for any operator space $A$, the space $\max(A)$
is completely isometric to a quotient space  of $\ell_1(I)$
for some set $I$ (\cf  [BP1]). 
Then since $\Cal K$ is nuclear, we have $d(\Cal K) = 2$, so that 
$\Cal K$ is completely isometric to a quotient 
space of $\Cal X_2(I)$ 
for some suitable countable set $I$ (see Remark 4.6). Therefore, $\Cal 
K\otimes_h \Cal K$ is completely isometric to a quotient space of 
$\Cal X_2(I) \otimes _h \Cal X_2(I) = \Cal X_4(I)$. Since $R$ and 
$C$ are completely isometric to quotients of $\Cal K$, it follows that 
$S_1 = R\otimes_h C$ is completely isometric to a quotient of $\Cal K 
\otimes_h \Cal K$. Finally, since every separable operator space is 
completely isometric to a quotient of $S_1$ (\cf  [B2, p. 24]), we conclude that 
every (resp.\ separable) operator space is completely isometric to a 
quotient of $\Cal X_4(I)$ for some set $I$ (resp.\ countable). The 
modification for the non-separable case is immediate.
\endremark

\head \S 6. $C^*$-algebras\endhead

\n Let $A$ be a unital $C^*$-algebra. The ``setting"
 used in this section is the same
as in \S 4, \ie $E=\max(A)$ and $\Cal A=A$. Note that in the
$C^*$-case, we have $\tilde A_1=A$ completely isometrically,
and (4.2) becomes $\Phi_A(c)=\Phi(i,c)$ for all $c>1$.

It is known that any nuclear $C^*$-algebra satisfies $d(A)\le 2$
(\cf  Bunce [Bu] and Christensen [C3]). In this section we study the
converse. Essentially, we show that if $A$ admits sufficiently many
type $II$ representations, then indeed the converse holds. More precisely
we will prove the following.

\proclaim{Theorem 6.1} Let $A$ be a $C^*$-algebra such that $d(A)\le
2$. Then, for any representation $\pi\colon \ A\to B(H)$ such that $\pi(A)$
generates   a semi-finite von~Neumann subalgebra of $B(H)$, the
bi-commutant $\pi(A)''$ is injective.
\endproclaim

\proclaim{Corollary 6.2} If $\dim H=\infty$, then
$d(B(H)) = 3.$
\endproclaim

\demo{Proof} Indeed, on one hand we know by Haagerup's result ([H1, Prop. 1.8]) 
that 
$d(B(H)) \le 3$.
On the other hand, by Joel Anderson's results in [A], there is a type
$II_\infty$ representation $\pi\colon\ B(H) \to B(\Cal H)$ such that
$\pi(B(H))'' \cong M\otimes B(H)$ where $M$ is a $II_1$ factor containing a
non trivial ultraproduct of matrix spaces. By 
Wassermann's result [W], we
know that the latter is not injective, so $\pi(B(H))''$ is not injective.
Thus, Theorem~6.1 implies in particular 
that $d(B(H)) \ge 3$. (I am grateful
to Simon Wassermann for kindly directing me  to 
 Anderson's result and  explaining to me 
its consequences.)\qed
\enddemo

\remark{Remark} By [C4], for any $II_1$-factor $M$ 
with property $\Gamma$ 
we have $ d(M)\le 44$. Since these cannot be 
nuclear ([W]),
the preceding result ensures that  
$3\le  d(M)\le 44$. It would
of course be interesting to reduce the interval of possible
values of $ d(M)$.
\endremark

The proof uses the following results.

The first lemma is a simple variant of a result from [JP].

\proclaim{Lemma 6.3} Let $A$ be any $C^*$-algebra. Then for any $n$
and any $\xi_1,\ldots, \xi_n \in A^*$ we have
$$\left(\sum \|\xi_i\|^2\right)^{1/2} \le 4\left\|\sum^n_{i=1} \xi_i
\otimes e_i\right\|_{A^*\otimes_{\text{\rm min}} \max(\ell_2)}. \tag6.1$$
\endproclaim

\demo{Proof} (The proof combines observations made independently by M.~Junge [J] 
and the
author.) Let $u\colon \ A\to \max(\ell^n_2)$ be the map defined by $u(a) =
\sum^n_1 \xi_i(a)e_i$. Let $E$ be a finite dimensional operator space. We
use the same notation as in [JP], i.e.\ we denote $d_{SK}(E) =
\inf\{\|v\|_{cb} \|v^{-1}\|_{cb}\}$ where the infimum runs over all
possible isomorphisms $v\colon\ E\to \widetilde E$ between $E$ and a
subspace of the   $C^*-$algebra 
 of all compact operators on $\ell_2$, 
which we have denoted above
by $\Cal K$. 

\n Let $a_1, \ldots, a_n$ be a finite subset of $A$ and let $E\subset A$ be
their linear span. Then the mapping $u_{|E}\colon \ E\to \max(\ell^n_2)$
factors through $A$ completely boundedly with a corresponding constant $\le
\|u\|_{cb}$. Fix $\epsilon>0$. By Lemma 6.2.11 in [P3] this implies that 
$u_{|E}$ can be written
as a composition $u_{|E} = u_2u_1$ with $u_1\colon \ E\to \widetilde E$ and
$u_2\colon \ \widetilde E\to \max \ell^n_2$ such that
$\|u_1\|_{cb} = 1$, $d_{SK}(\widetilde E) = 1$ and $\|u_2\|_{cb} \le
\|u\|_{cb}(1+\epsilon)$. By the main result in [JP], this implies
$$\align
\left|\sum\nolimits^n_1 \xi_i(a_i)\right| &=
\left|\sum\nolimits^n_1\langle u(a_i),e_i\rangle\right|\\
&= \left|\sum\nolimits^n_1 \langle u_2u_1(a_i),e_i\rangle\right|\\
&\le 4~ d_{SK}(\widetilde E) \|u_2\|_{cb} \left(\sum
\|u_1(a_i)\|^2\right)^{1/2}\\
&\le 4\left(\sum \|a_i\|^2\right)^{1/2} \|u\|_{cb}(1+\epsilon).
\endalign$$
Hence, since  $\epsilon>0$ is arbitrary, and since
$$  \|u\|_{cb} = \left\|\sum^n_{i=1} \xi_i
\otimes e_i\right\|_{A^*\otimes_{\text{\rm min}} \max(\ell_2)},$$
  taking the supremum over all possible $n$-tuples
$(a_i)_{i\le n}$ in
$A$, we obtain (6.1).\qed
\enddemo

\proclaim{Lemma 6.4} Let $(e_i)$ be the canonical basis of the operator
space $\max(\ell_2)$. Let $H$ be any Hilbert space and let $X$ be either
$B({\comp},H)$ or $B({ H}^*,{\comp })$, or equivalently let $X$ be
either the column Hilbert space or the row Hilbert space. Then for all
$x_1,\ldots, x_n$ in $X$ we have
$$\left\|\sum\nolimits^n_1 x_i\otimes e_i\right\|_{X\otimes_{\text{\rm min}}
\max(\ell_2)} \le \left(\sum \|x_i\|^2\right)^{1/2}.$$
\endproclaim

\demo{Proof} Assume $X = B({\comp },H)$ or $B(H^*,{\comp })$. We identify $X$ 
with
$H$ as a vector space.	Let
$(\delta_m)$ be an orthonormal basis in $H$. Observe that for any finite
sequence $a_m$ in $B(\ell_2)$ we have in both cases
$$\left\|\sum \delta_m \otimes a_m\right\|_{\text{\rm min}}\le \left(\sum
\|a_m\|^2\right)^{1/2}.$$
whence we have, for any $x_1,\ldots, x_n$ in $X$,
$$\align
\left\|\sum x_i\otimes e_i\right\| &= \left\|\sum_m \delta_m
\otimes \sum_i \langle x_i,\delta_m\rangle e_i\right\|\\
&\le \left(\sum_m \left\|\sum_i \langle x_i,\delta_m\rangle
e_i\right\|^2\right)^{1/2}\\
&= \left(\sum_{m, i} |\langle x_i,\delta_m\rangle|^2\right)^{1/2} =
\left(\sum_i \|x_i\|^2\right)^{1/2}.
\endalign$$\qed
\enddemo

\demo{Proof of Theorem 6.1} Recall that, since $\pi(A)$ is a quotient
$C^*$-algebra of $A$, we have obviously $d(\pi(A)) \le d(A)$. Hence
it suffices to prove the statement with $\pi(A)$ in the place of $A$. More
precisely, we assume given $A\subset B(H)$ such that $M=A''$ admits a
faithful semi-finite normal trace denoted by $\tau$ and we must show that
$d(A) \le 2$ implies that $M$ is injective. First we can reduce to the
finite case:\ indeed it suffices to show that, for any projection $p$ in
$M$ and $0 < \tau(p)<\infty$, the algebra $pMp$ is injective. Then, by a
result due to Connes for factors and to Haagerup [H2] in the general case,
$pMp$ is injective iff there is a constant $C$ such that for any central
projection $q\not= 0$ in $pMp$, for any $n$ and any $n$-tuple $u_1,\ldots, u_n$ 
of
unitaries in  $pMp$, we have
$$n \le C\left\|\sum\nolimits^n_1 (qu_i) \otimes
\overline{(qu_i)}\right\|_{\text{\rm min}}.\tag6.2$$
Fix  $p,q$ and $u_1,\ldots, u_n$ unitary in $pMp$ as above. We will show
that (6.2) holds.
\enddemo

\n Let $\xi_i\in A^*$ be the functional defined by 
$$\xi_i(a) = \tau(qu_ia)$$
and let
$$t = \left\|\sum\nolimits^n_1 (qu_i) \otimes
\overline{(qu_i)}\right\|_{\text{\rm min}}^{1/2}. $$
Let $\sigma\colon \ M\to B(L^2(M,\tau))$ be the classical representation of
$M$ as left multiplications on $L^2(M,\tau)$. Clearly $\left\|\sum
\sigma(qu_i) \otimes \overline{\sigma(qu_i)}\right\|_{\text{\rm min}} \le
\left\|\sum qu_i \otimes \overline{qu_i}\right\|_{\text{\rm min}}$ hence ([P6, 
Cor. 
2.7]) there is
a decomposition
$$\sigma(qu_i) = a_i+b_i\tag6.3$$
in $B(L^2(M,\tau))$ with
$$\sum a^*_i a_i \le t^2\quad \text{and}\quad \sum b_ib^*_i \le
t^2.$$ Consider the mapping $u\colon \ A\to \max(\ell_2^n)$ defined by
$$\forall a\in A\quad u(a)=\sum_1^n \xi_i(a) e_i.$$
Then (using the tracial property of $\tau$) we have (recall $qu_i = qu_iq
=u_iq$)
$$\align
u(ab) &= \sum\nolimits^n_1 \tau(bqu_ia) e_i\\
&= \sum\nolimits^n_1 \tau(bq(u_i)qa) e_i.
\endalign$$
Let us denote by $j\colon \ A\to L^2(M,\tau)$ the mapping defined by $j(a)
= qa$. We can write for all $a,b$ in $A$
$$u(ab) = \sum\nolimits^n_1 \langle(a_i+b_i)j(a), j(b^*)\rangle e_i$$
(here the scalar product is in $L^2(M)$) hence
$$u(ab) = \varphi_1(a,b) + \varphi_2(a,b)$$
where
$$\varphi_1(a,b) = \sum\nolimits^n_1 \langle a_ij(a), j(b^*)\rangle
e_i$$ and
$$\varphi_2(a,b) = \sum\nolimits^n_1 \langle b_ij(a), j(b^*)\rangle
e_i.$$
Now we assume $\max(\ell^n_2) \subset B(\Cal H)$ completely isometrically
and we claim that for some Hilbert space $\widehat H$ there are mappings
$$\beta_i\colon \ A\to B(\Cal H,\widehat H)\quad \text{and}\quad
\alpha_i\colon \ A\to B(\widehat H,\Cal H)$$
for $i=1,2$ such that
$$\varphi_i(a,b) = \alpha_i(a) \beta_i(b)$$
and
$$\|\alpha_1\| \, \|\beta_1\| \le  \tau(q)t,\quad \|\alpha_2\|\, \|\beta_2\|
\le \tau(q)t.$$
Taking this for granted for the moment, let us now complete the argument.
By Theorem~4.2  our assumption $d(A)\le 2$ implies that
there is a constant $K$ such that the product mapping $P\colon \ \max(A)
\otimes_h \max(A) \to A$ satisfies for all maps $u\colon \ A\to B(\Cal
H)$
$$\|u\|_{cb} \le K\|uP\|_{cb}.$$
Note that if $\hat u\colon \ \max(A) \times \max(A) \to B(\Cal H)$ is the
bilinear form associated to $uP$, then by [CS1-2, PaS] we know that
$$\|uP\|_{cb}= \|\hat u\|_{cb}$$
and since $\hat u = \varphi_1 +\varphi_2$ we obtain
$$\|u\|_{cb} \le K[\|\hat u\|_{cb}] \le K(\|\varphi_1\|_{cb} +
\|\varphi_2\|_{cb})$$
but by the specific factorization of $\varphi_1$ and $\varphi_2$ given
above we have $\|\varphi_i\|_{cb} \le \|\alpha_i\|\, \|\beta_i\|$ whence
$$\align
\|u\|_{cb} &= K(\|\alpha_1\|\, \|\beta_1\| + \|\alpha_2\|\,
\|\beta_2\|)\\
&\le 2Kt \tau(q).
\endalign$$
Equivalently we have
$$\left\|\sum \xi_i\otimes e_i\right\|_{A^* \otimes_{\text{\rm min}}
\max(\ell^n_2)} \le 2\tau(q)tK.$$
By (6.1) this implies
$$\left(\sum \|\xi_i\|^2_{A^*}\right)^{1/2} \le 8\tau(q) tK.$$
But on the other hand
$$\|\xi_i\|_{A^*} \ge |\xi_i(u^*_i)| = \tau(qu_iu^*_i) = \tau(q)$$
hence
$$\tau(q) \sqrt n \le 8\tau(q )tK$$
and finally $n\le 64 K^2 t^2$. Thus we obtain (6.2).

This completes the proof modulo  the claim. We now turn to the latter
claim. Let $L = L^2(M,\tau)$. We denote by $x\to r(x)\in B(\bar L,{\comp})$
 the canonical identification. Note that $r(x) r(y)^* \in B({\comp},{\comp })$ 
can be identified with $\langle x,y\rangle$.

\n With this identification, we have for all $a,b$ in $A$
$$\align
\varphi_1(a,b) &= \sum\nolimits^n_1 r(a_ij(a)) r(j(b^*))^*
e_i\\
&= \left(\sum\nolimits^n_1 r(a_ij(a)) \otimes e_i\right) \circ
(r(j(b^*))^* \otimes I).
\endalign$$
We set $\alpha_1(a) =  \sum^n_1 r(a_ij(a)) \otimes e_i$ and
$\beta_1(b) = r(j(b^*))^* \otimes I.$
Then $\|\beta_1(b)\| = \|j(b^*)\|_L = \|qb^*\|_{L^2(\tau)} \le
\|b\| \tau(q)^{1/2}$ and by Lemma~6.4
$$\align
\|\alpha_1(a)\| &\le \left(\sum\nolimits^n_1
\|r(a_ij(a))\|^2_L\right)^{1/2}\\
&\le \left(\sum\nolimits^n_1 \|a_ij(a)\|^2_L \right)^{1/2}\\
&\le \left\|\sum a^*_ia_i\right\|^{1/2} \|j(a)\|_L\\
&\le t\tau(q)^{1/2} \|a\|.
\endalign$$
Hence we obtain $\|\alpha_1\|\, \|\beta_1\|\le t\tau(q)$ as announced.
Similarly we define
$$\alpha_2(a) = r(j(a)) \otimes I$$
and
$$\beta_2(b) = \sum\nolimits^n_1 r(b^*_ij(b^*))^* \otimes e_i.$$
Then clearly $\varphi_2(a,b) = \alpha_2(a) \beta_2(b)$ and this time we
have $\|\alpha_2\| \le \tau(q)^{1/2}$ and
$\|\beta_2\| \le t\tau(q)^{1/2}$, whence $\|\alpha_2\|\, \|\beta_2\|\le
t\tau(q)$.
This completes the proof of the claim, and also of Theorem~6.1.\qed

In particular, since
$C^*(G)$  or  $C_\lambda^*(G)$ is nuclear iff 
$G$ (discrete) is amenable, (\cf  e.g. [La])
we recover some results from \S 3, as follows.

\proclaim{Corollary 6.5} Let $G$ be a discrete group and let $A$
be either $C^*(G)$ or the 
reduced $C^*$-algebra $C_\lambda^*(G)$. Then $d(A)\le 2$ iff $G$ is amenable.
\endproclaim

\remark{Remark} Note however, that the equivalence with (v) in Theorem 3.2
concerning the spaces of coefficients 
does not follow from this new approach.
\endremark

\proclaim{Corollary 6.6} Let $A$ be a   
$C^*$-algebra which generates
a non-injective semi-finite von Neumann algebra. Then for any
$c>1$, there is a unital homomorphism
$u_c\colon\ A\to B(H)$ with $\|u_c\|\le c$ and  $\|u_c\|_{cb}\ge c^3 $.
\endproclaim

\remark{Remark 6.7} A unital  $C^*$-algebra $A$ satisfies the similarity
property (\ie $d(A)<\infty$) as soon as  $\Phi_A(c) <\infty$ for {\it some}
$c>1$. Indeed, this follows from Lemma 2.3 and the remark preceding it.
\endremark

\remark{Remark} The following result proved in [H1] and [C3] plays an
important r\^ole in these papers:\ Let $u\colon\ A_1\to A_2$ be a bounded
homomorphism between $C^*$-algebras. Then for any finite subset
$(x_i)$ in $A_1$ we have
$$\left\|\sum u(x_i)^* u(x_i)\right\|^{1/2} \le \|u\|^2 \left\|\sum
x^*_ix_i\right\|^{1/2}.$$
The next result shows that the exponent 2 cannot be improved in this result.
\endremark

\proclaim{Proposition 6.8} Suppose that a number $\alpha\ge 1$ has the
following property:\ there is a constant $K$ such that for any bounded
homomorphism
$$u\colon \ B(H)\to B(H)\qquad (\dim H=\infty)$$
and for any finite subset $x_1,\ldots, x_n$ in $B(H)$ we have
$$\left\|\sum u(x_i)^* u(x_i)\right\|^{1/2} \le K\|u\|^\alpha
\left\|\sum x^*_i x_i\right\|^{1/2}.$$
Then necessarily $\alpha\ge 2$.
\endproclaim

\demo{Proof} Our assumption can be written as follows:\ for any $c\ge 1$
and any unital homomorphism $u$ with $\|u\|\le c$,
 we have for any $n$ and any 
$x_1,\ldots, x_n$ in $B(H)$
$$\left\|(I_{\Cal K}\otimes u) \left(\sum e_{i1}\otimes
x_i\right)\right\|_{\text{\rm min}}\le Kc^\alpha \left\|\sum
e_{i1} \otimes x_i\right\|_{\text{\rm min}}.$$ 
In other words, the subspace $X$ spanned in $\Cal K$
by the sequence $(e_{i1})$ $(i=1,2,...)$ satisfies
the assumption (2.11) in Theorem 2.6. Assume $\alpha<2$.
Then, by
Theorem~2.6,  (2.11) actually holds for $\alpha=1$ (for some $K$).
Thus, if $\alpha<2$ we may as well assume $\alpha=1$. But then
Haagerup's argument in [H1] (or the proof presented in [P1,
chapter~7]) will lead to $d(B(H)) \le 2$, which
contradicts Corollary~6.2. Thus we must have $\alpha\ge
2$.\qed
\enddemo

\head \S 7. The Blecher-Paulsen factorization\endhead

\n  In this section, we connect our description of the enveloping algebra
$\tilde A_1$
with some ideas of Blecher and Paulsen in [BP2]. We take
 a slightly more general viewpoint than them in 
order to cover the
situation of a group (or an algebra) generated by a subset, but
the main idea is in [BP2].

We consider our usual ``setting" $(i,E, \Cal A)$, where
$E$ is an operator space, $\Cal A$ a unital operator 
algebra (not assumed complete) , and
$i\colon\ E\to \Cal A$ is a completely contractive
 linear injection with range
generating $\Cal A$. But in addition we will
 assume throughout
this section that $E$ is ``unital", by which we mean that $E$ contains a norm
one element $e$ such that $i(e) = 1_{\Cal A}$.
 
Consider again the
algebra $\tilde A_1$ as defined above, with unital embeddings $E\subset
\Cal A\subset \tilde A_1$.

\n It will be convenient to consider $E$ as ``included"
into $\Cal A$ and to view $i$ as an inclusion map.
 The reader
should be warned however that $i$ will
  generally not be assumed completely
isometric: in general the operator space structure
on $\Cal A$ only plays  a auxiliary r\^ole. 
What really matters here
is the given operator space structure on $E$ 
and the resulting operator 
algebra one on $\tilde A_1$, which appears as ``generated"
 by $E$.

\proclaim{Theorem 7.1}
With the above notation, let $n$ be a positive integer. 
Then the
following properties of
an element $x$ in $M_n(\Cal A)$ are equivalent:
\roster
  \item"{\text{\rm (i)}}" $\|x\|_{M_n(\tilde A_1)} <1$.
  \item"{\text{\rm (ii)}}" The matrix $x$ can be written, for some
integer $N$ and some integer $d$, as a matricial product of
the form
$$x = \alpha_0 D_1\alpha_2D_2\ldots D_d\alpha_{d}$$
where $\alpha_0\in M_{n N}$, $\alpha_1\in M_{N}$,..., $\alpha_{d-1}\in M_{N}$,
$\alpha_{d}\in M_{Nn}$ are scalar matrices (\ie $\alpha_0$ and $\alpha_{d}$
are  rectangular of size $n\times N$ and $N\times n$, and the others are square
matrices of size $N\times N$),
 and
$D_1,\ldots, D_d$ are  $N\times N$  matrices with entries in $E$,
and finally we have
$$\prod^{d}_{i=0} \|\alpha_i\|  \prod^d_{i=1} \|D_i\|_{M_N(E)} <
1.$$
\endroster
\endproclaim

\demo{Proof} The proof follows from an immediate adaptation of
an argument in [BP2]. We merely sketch it.
It is clear that (ii) implies (i).
 Conversely assume (i). This means that
there is a number $\theta<1$, such that for any contractive
unital homomorphism
$u\:  \Cal A \to B(H)$,
 we have $$\|[I_{M_n}\otimes u](x)\|_{M_n(B(H))}\le \theta.\tag7.1$$
Now, for any $n$ and for any $x$ in $M_n(\Cal A)$, let us denote
by $\|x\|_{(n)}$ the infimum of $$\prod^{d}_{i=0} \|\alpha_i\|  \prod^d_{i=1} 
\|D_i\|$$
over all possible factorizations (with arbitrary $N$ and $d$) of $x$ as in (ii) 
above. 
Then by the main result in [BRS], there is a unital
homomorphism $u\:  \Cal A \to B(\Cal H)$ 
with the help of which
$\|.\|_{(n)}$ can be identified
 with $\|u(.)\|_{M_n(B(\Cal H))}$ for all $n$. 
But then, clearly $ui$ is completely contractive,
 since we trivially have (by definition of $\|.\|_{(N)}$)
$\|[I_{M_N}\otimes u](D)\|_{M_N(B(H))}=
\|D\|_{(N)}\le \|D\|_{M_N(E)}$ 
for any $D$ in $M_N(E)$ and any $N$.
Hence, returning
to the particular $n$ and $x$ 
appearing in (i), by (7.1) we must have 
 $\|x\|_{(n)}\le \theta$. Equivalently,
we obtain (ii). \qed
\enddemo

\remark{Remark 7.2} Recall that when $A$ is a $C^*$-algebra 
(resp. $A=C^*(G)$) and 
$E=\max(A)$ (resp. $E=\ell_1(G)$) as in \S 5 (resp. \S 3), then $\tilde A_1=A$ 
completely
isometrically.   When the latter holds,
 Theorem 7.1 gives a
characterization of the 
elements of the unit ball of $M_n(A)$.
 Note also that when $E=\max(E)$, the elements of $M_N(E)$
admit a specific factorization (a kind of diagonalization) 
described above in (0.4).
\endremark

As  application, we have the following apparently
new characterization  of  the coefficients of unitary representations of
a group $G$, \ie of the elements of the space
$B(G)$, as follows. (Take $H$ unidimensional in the next
statement, then  (i) below is the same as
saying that the norm of
$f$ in $B(G)$ is $\le K$.)

\proclaim{Corollary 7.3} Let $G$ be any discrete group,
and let $\Gamma\subset G$ be a  subset
containing the unit element and 
 generating $G$ in the sense
that every element of $G$ can be
 written as a product of elements of
$\Gamma$. Let $K\ge 0$ be a fixed
constant. The following properties of a function $f\colon
\ G\to B(H)$ are equivalent:
\roster 
\item"{\text{\rm (i)}}" There are a
unitary representation $\pi\colon \ G\to B(H_\pi)$ and
operators $\xi\colon \ H_\pi \to H$ and $\eta\colon\ H\to
H_\pi$ such that $f(t) = \xi\pi(t)\eta$ for any $t$ in $G$
and $\|\xi\|\, \|\eta\| \le K.$ 
\item"{\text{\rm (ii)}}" For each
$N\ge 1$, the function $f_N\colon \ \Gamma^N\to B(H)$
defined by $f_N(t_1,\ldots, t_N) = f(t_1t_2\ldots t_N)$
extends (with the obvious identification) 
to an element of \break
$cb(\ell_1(\Gamma) \otimes_h \cdots \otimes_h
\ell_1(\Gamma)$, $B(H))$ (where the tensor product is
$N$-fold) with norm $\le K$.
\item"{\text{\rm (iii)}}" Same as
(ii) with $\Gamma=G$.
\item"{\text{\rm (iv)}}" For any $N\ge 1$,
there are bounded functions $F_i\colon \ \Gamma \to
B(H_{i+1}, H_i)$, with  $\sup_\Gamma\|{F_i}\|\le
1$ for all $i$,  where $H_i$ are Hilbert spaces with
$H_{N+1}=H$ and $H_1=H$, such that
$$f(t_1t_2\ldots t_N) =K\   F_1(t_1) F_2(t_2) \ldots
F_N(t_N).\leqno \forall\ t_1,\ldots, t_N\in \Gamma$$
\endroster
\endproclaim

\demo{Proof} By the factorization of c.b.\ maps (\cf  [Pa1]) (i) holds iff $f$
extends linearly to a c.b.\ mapping $u\colon \ C^*(G)\to B(H)$ with
$\|u\|_{cb}\le K$. We consider now the setting $(i,E,\Cal A)$ defined by
$E = \ell_1(\Gamma)$, $\Cal A = \ell_1(G)$ (viewed as a subalgebra of
$C^*(G)$). Then, as already mentioned in Remark 7.2, for any $n$, we have an 
isometric identity
 $M_n(\tilde A_1)
= M_n(C^*(G))$. 
Let us denote as in \S 3, by $W_N$ the natural
product map from $E_N =\ell_1(\Gamma)\otimes_h...\otimes_h \ell_1(\Gamma)$ ($N$ 
times)
into $C^*(G)$. Then, by Theorem 7.1, we have $\|u\|_{cb}\le K$ iff 
$$\sup_N \|u W_N\|_{cb(E_N,B(H))} \le K.\tag7.2$$
But now,  (7.2) is but a reformulation
of (ii), so that (i) is equivalent to (ii). Moreover, since 
(i)~$\Rightarrow$~(ii) is valid
for any $\Gamma$, it holds when $\Gamma=G$, whence
(ii)~$\Rightarrow$~(iii), and the converse is obvious. Finally, the equivalence
between (ii) and (iv) follows from the well known 
factorization theorem of $c.b.$ multilinear maps (\cf  [CS1-2,
PaS]). \qed 
\enddemo

\head \S 8.  Banach algebras\endhead

\n The general method of this paper can be applied in   other situations when 
studying
a Banach algebra $A$ given together with a generating system, or a family of 
generating
subalgebras. The r\^ole of the ``degree" is then played by the minimal length of 
the products 
necessary to generate (in a suitable Banach algebraic sense) the unit  ball (or 
some ball
centered at the origin).  One can also develop our approach for a   general 
``variety of Banach
algebras" (in the sense of [Dix])
instead of that of operator algebras. To illustrate briefly what we have in 
mind,   take the
variety
of all Banach algebras, then our basic idea leads to:

\proclaim{Theorem 8.1}  Let $A$ be a Banach algebra with unit ball $B_A$. 
Consider a  subset
$\beta\subset B_A$ and assume that the algebra it generates, denoted by $\Cal
A$, is dense in
$A$. For $n=1,2,...$, let us denote by $\beta^n$ the set of all products of $n$ 
elements taken in
$\beta$. Let $d$ be a positive  integer and let $\alpha$ be any number such that 
$d\le \alpha
<d+1$. Consider the following properties:
\roster
\item"{${\text{\rm (i)}}_\alpha$}" There is a constant $K$ such that, for any  
Banach 
algebra $B$ and  for
any homomorphism
$u\: \Cal A\to B$, if $u$ is bounded on $\beta$, $u$ is continuous and we have
$$\|u\|\le K\sup_{x\in \beta} \|u(x)\|^\alpha.$$ 
\item"{${\text{\rm (ii)}_d}$}" There 
is a constant $K'$ such that (here ${\text{\rm
aconv}}$ stands
for the absolutely convex hull)
$$B_A\subset K' \overline{\text{\rm
aconv}}(\beta\cup 
\beta^2\cup...\cup\beta^d).$$
\item"{${\text{\rm (iii)}}_d$}" There is a constant $K''$ such that, for any  
Banach 
algebra $B$ and  for
any continuous homomorphism
$u\: \Cal A\to B$, we have
$$\|u\|\le K''\sup_{x\in \beta} \|u(x)\|^d.$$ 
Then ${\text{\rm (i)}_\alpha}\Rightarrow {\text{\rm (ii)}_d}$. 
Moreover, if $\beta$ contains a unit element for
$A$, then conversely ${\text{\rm (ii)}_d}\Rightarrow {\text{\rm (iii)}_d}$.
\endroster
\endproclaim

\demo{Proof} (Sketch) In this proof, we will say ``morphism" for homomorphism 
with values 
in a Banach
algebra. We will follow the same strategy as in \S 1 and \S 2. Let 
$\|u\|_\beta=\sup\{ \| u(x)\| \
|
\ x\in
\beta\}$. Let $c\ge 1$ and let $\Cal C_c$ be the set of all morphisms
$u\: \Cal A\to B_u$ with  $\|u\|_\beta \le c$. We define the Banach algebra 
$\tilde A_c$
as the completion of $\Cal A$ for the embedding $J\: \Cal A\to \oplus_{u\in 
\Cal C_c} B_u$
defined by $J(x)=\oplus_{u\in \Cal C_c} u(x)$.
Let $F_\beta$ be the free semi-group with free generators indexed by $\beta$. We 
will consider $\beta$
as a subset of $F_\beta$. Consider the space $\ell_1(F_\beta)$, viewed as a 
Banach algebra for
convolution. Let $\delta_t$ ($t\in F_\beta$) denote the  canonical basis and let 
$\Cal B$
be the dense subalgebra linearly generated by  $\delta_t$ ($t\in F_\beta$), 
equipped with the
induced norm. 

\n Let $c\ge 1$ and $z=1/c$. We have a unique morphism $\pi_z\: \Cal B\to 
\Cal A$
such that $\pi_z(\delta_x)= z x$ for all $x$ in $\beta$.
It is easy to check that if 
$u\: \Cal A\to B_u$ is any morphism, then $\|u\|_\beta \le c$ iff $\|u\pi_z \| 
\le 1$.
Moreover, $\tilde A_c$ can be identified with the completion of $\Cal 
B/\ker(\pi_z)$,
and $ (i)_\alpha$ means that 
$$\forall x \in \Cal A\quad  \|x\|_{\tilde A_c}\le Kc^\alpha \| x\|_A.$$
 The proof can then be completed by arguing as in Theorem 2.5. We leave the 
remaining details to the reader.
\qed
\enddemo
\medskip

\n {\bf Acknowledgement:} I am very grateful to
 Christian Le Merdy
for pointing out many defects of the preliminary versions.


\bigskip


\Refs
\widestnumber\key{MPSZ}
 
\ref \key  A \by J.  Anderson \paper Extreme points
 in sets of linear functionals in $\Cal B(\Cal H)$
\jour J. Funct. Anal. \vol 31 \yr1979 \pages 195--217\endref

\ref \key B1 \by D. Blecher \paper A completely bounded
characterization of operator algebras \jour Math. Ann.  
\vol303 \yr1995 \pages 227--240\endref

\ref \key B2 \bysame \paper The standard
 dual of an operator space
 \jour Pacific J. Math. \vol 153 \yr1992 \pages 15--30\endref


\ref \key BRS \by D. Blecher, Z.J. Ruan  and  A. Sinclair \paper 
A characterization of operator algebras \jour J. Funct. Anal.
\vol 89 \yr1990  \pages 188--201\endref


\ref \key BP1 \by D. Blecher and V. Paulsen \paper  Tensor products of
operator spaces  \jour J. Funct. Anal.  \vol 99 \yr1991 \pages 262--292\endref

\ref \key BP2 \bysame \paper   Explicit
construction of universal operator algebras and
applications to polynomial factorization \jour Proc. Amer.
Math. Soc. \vol  112 \yr1991 \pages 839--850\endref 

\ref \key BS \by D. Blecher  and R. Smith \paper The dual of the Haagerup tensor
 product \jour Journal London Math. Soc. \vol 45 \yr1992 \pages 126--144\endref
 

\ref \key Bou  \by J. Bourgain \paper  On the similarity problem for
polynomialy bounded operators on Hilbert space \jour Israel J.
Math. \vol 54 \yr1986 \pages 227--241\endref

\ref \key Bo1 \by  M. Bo$\dot{z}$ejko \paper  Remarks on Herz-Schur
multipliers on free groups \jour  Mat. Ann. \vol 258 \yr1981
\pages 11--15\endref


\ref \key Bo2  \bysame \paper  Positive definite bounded
matrices and a characterization of amenable groups \jour 
Proc. A.M.S. \vol  95 \yr1985 \pages 357--360\endref 


\ref \key BF1 \by   M. Bo$\dot{z}$ejko and  G. Fendler \paper  Herz-Schur
multipliers and completely bounded multipliers of the
Fourier algebra of a locally compact group \jour  Boll. Unione
Mat. Ital.  \vol (6) 3-A \yr1984 \pages 297--302\endref


\ref \key BF2 \bysame \paper Herz-Schur
multipliers and uniformly bounded representations of
discrete groups \jour Arch. Math. \vol 57 \yr1991 \pages 290--298\endref 

\ref \key Bu  \by J. W. Bunce \paper  The similarity problem for
representations of $C^*$-algebras \jour  Proc. Amer. Math.
Soc. \vol  81 \yr1981 \pages 409--414\endref

\ref \key CE \by M.D.  Choi and E. Effros \paper  Nuclear
C*-algebras and injectivity: The
       general case \jour Indiana Univ. Math. J. \vol 26 \yr1977
\pages 443--446\endref

\ref \key C1  \by   E. Christensen \paper   Extensions of derivations \jour 	
J.
Funct. Anal. \vol 27 \yr1978 \pages 234--247\endref

\ref \key C2  \bysame \paper  Extensions of derivations II \jour
Math. Scand. \vol 50 \yr1982 \pages  111--122\endref

 

\ref \key C3  \bysame \paper On non self adjoint
representations of operator algebras \jour Amer. J. Math. \vol
103 \yr1981 \pages 817--834\endref

\ref \key C4 \bysame   \paper Similarities of $II_{1}$
factors with property $\Gamma$ \jour  Journal Operator Theory 
\vol 15 \yr1986 \pages 281--288\endref 
 
\ref \key CES \by  E. Christensen,  E. Effros  and  A. Sinclair \paper 
Completely bounded multilinear maps and
$C^*$-algebraic cohomology \jour Invent. Math. \vol 90 \yr1987
\pages 279--296\endref

\ref \key CS1 \by E. Christensen and A. Sinclair \paper
Representations of completely bounded multilinear operators \jour
J. Funct. Anal. \vol 72 \yr1987 \pages 151--181\endref


\ref \key CS2  \bysame  \paper A survey of
completely bounded operators \jour  Bull. London Math. Soc. \vol 
21 \yr1989 \pages 417--448\endref

 

\ref \key Co \by   M. Cowling \paper  Uniformly bounded representations
of Lie groups \jour CIME Course Lecture Notes\endref

 \ref \key CL \by P.C. Curtis and R.J. Loy \paper A note on 
amenable algebras of operators \jour
Bull. Australian Math. Soc.  \vol 52 \yr1995 \pages 327--329\endref
 
 \ref \key Da \by  A.M. Davie \paper  Quotient algebras of uniform
algebras \jour
J. London Math. Soc.  \vol 7 \yr1973 \pages 31--40\endref

 
\ref \key DCH  \by J. de Canni\`ere  and U.  Haagerup \paper  Multipliers
of the Fourier algebras of some simple Lie groups and
their discrete subgroups \jour Amer. J. Math. \vol  107 \yr1985 \pages
455--500\endref

\ref \key Di \by J. Dixmier \paper  Les moyennes invariantes dans les
semi-groupes et leurs applications \jour  Acta Sci. Math.
Szeged  \vol 12 \yr1950 \pages 213--227\endref

 \ref\key Dix \by P. Dixon \paper Varieties of
 Banach algebras \jour
 Quarterly J. Math. Oxford \vol 27 \yr1976 \pages 481--487\endref

 

  \ref \key E \by E. Effros \paper Amenability and virtual diagonals 
for von Neumann algebras \jour
J. Funct. Anal. \vol 78 \yr1988 \pages 137--153\endref
 

 \ref \key EM \by
 L. Ehrenpreis  and F.I. Mautner \paper  Uniformly bounded
representations of groups \jour Proc. Nat. Acad. Sc. U.S.A. \vol 41
\yr1955 \pages 231--233\endref

 

 
 
 \ref \key ER  \by E. Effros     and Z.J. Ruan \paper  A new approach to
operators spaces \jour
 Canadian Math. Bull. \vol 
34 \yr1991 \pages 329--337\endref


\ref \key Ey \by P. Eymard \paper    L'alg\`ebre de Fourier d'un groupe
localement compact \jour  Bull. Soc. Math. France  \vol 92 \yr1964 \pages
181--236\endref

 
 
 
 \ref \key FTP \by A. Fig\`a-Talamanca  and M.
Picardello \book   Harmonic Analysis on Free groups
\publ Marcel Dekker \publaddr New-York \yr 1983\endref

 


\ref \key Fo \by S. Foguel \paper  A counterexample to a problem of
Sz. Nagy. \jour Proc. Amer. Math. Soc. \vol 15 \yr1964 \pages 788--790\endref

\ref \key Ga  \by D.J.H. Garling \paper
 On the dual of a proper uniform algebra \jour
Bull. London Math. Soc. \vol 21 \yr1989 \pages 279--284\endref 


 

 

 
\ref \key H1 \by U. Haagerup \paper  Solution of the similarity
problem for cyclic representations of $C^*$-algebras \jour 
Annals of Math. \vol 118 \yr1983 \pages 215--240\endref


 \ref \key H2  \bysame \book Injectivity and decomposition of completely
 bounded maps in ``Operator algebras and their connection with Topology and
 Ergodic Theory'' \publ Springer Lecture Notes in Math. 1132
\yr1985 \pages 170--222\endref
 
\ref\key H3 \bysame \paper $M_{0} A(G)$ functions which are
not coefficients of
uniformly bounded representations \paperinfo Handwritten manuscript 
    \yr1985\endref


\ref \key H4  \bysame \paper 
 All  nuclear $C^*$-algebras are amenable \jour
 Invent. Math. \vol 74 \yr1983 \pages 305--319\endref



 \ref \key J \by M. Junge \book Factorization theory
 for spaces of operators \publ
Habilitationsschrift \publaddr Kiel \yr1996\endref

 
 
 \ref \key JP \by M. Junge and G. Pisier \paper Bilinear
 forms on exact operator
spaces and $B(H)\otimes B(H)$ \jour Geometric and 
Functional Analysis (GAFA Journal)  \vol
5 \yr1995 \pages 329--363\endref 

\ref \key Ka \by R. Kadison \paper  On the orthogonalization of
operator representations \jour Amer. J. Math. \vol 77 \yr1955 \pages 600--620
\endref

 
\ref \key KaR \by R. Kadison  and J. Ringrose \book Fundamentals of
the theory
of operator algebras \vol II  \publ Academic Press \publaddr 
New-York  \yr1986\endref

\ref \key Ki \by E. Kirchberg \paper  The derivation and the similarity problem 
are equivalent \jour
 J. Operator Th. \vol 36 \yr1996 \pages 59--62\endref
 
\ref\key Kis \by S.  Kisliakov \paper The proper uniform algebras are 
uncomplemented 
\jour Lomi preprint  \yr1988\endref
 


\ref\key KS   \by R.A. Kunze  and  E. Stein \paper  Uniformly bounded
representations and Harmonic Analysis of the $2\times
2$ real unimodular group \jour  Amer. J. Math.  \vol 82 \yr1960
\pages 1--62\endref

\ref \key La \by C. Lance  \paper  Tensor products and nuclear $C^*$-algebras.
Operator algebras and applications \jour Amer. Math. Soc.
Proc. Symposia Pure Math
  \vol 38 part 1 \yr1982 \pages 379--399\endref

 

\ref \key Le \by A. Lebow \paper A power bounded operator which is not
polynomially bounded \jour  Mich. Math. J. \vol 15
\yr1968\pages 397--399\endref
   

 
 \ref \key MPSZ\by A.M. Mantero,  T. Pytlik, R. Szwarc 
and A. Zappa \paper  Equivalence of two series of spherical representations of a 
free group \jour
Ann. di Matematica pura ed applicata (IV) \vol CLXV \yr1993 \pages 23--28\endref

 
 \ref \key O1\by  A. Yu. Ol'shanskii \paper
 On the problem of the existence of an invariant mean
 on a group \jour
Russian Math. Surveys \vol 35 \yr1980 \pages \pages 180--181\endref

 \ref \key O2\bysame \book
Geometry of defining relations in groups
\bookinfo \publ Kluwer\publaddr
Dordrecht, Netherlands\yr1991\endref

\ref \key Pat\by   A. Paterson \paper  Amenability
\jour  Amer. Math. Soc. Math. Surveys \vol 29
\yr1988\endref

\ref \key Pa1 \by V. Paulsen \book   Completely bounded maps and
dilations \bookinfo  Pitman Research Notes in
Math. \vol 146 \publ Longman, Wiley 
\publaddr New York  \yr1986\endref

\ref \key Pa2  \bysame \paper Completely bounded
maps on
$C^*$-algebras and invariant operator ranges \jour  Proc.
Amer. Math. Soc. \vol 86 \yr1982 \pages 91--96\endref 

\ref \key Pa3  \bysame \paper  Every completely 
polynomially bounded operator is similar to a
contraction \jour J. Funct. Anal. \vol 55 \yr1984 \pages 1-17\endref

\ref \key Pa4 \bysame \paper  Completely bounded homomorphisms
of operator algebras \jour  Proc.
Amer. Math. Soc. \vol 92 \yr1984 \pages 225--228\endref
 
\ref \key Pa5 \bysame \paper Representation of
Function algebras, Abstract
 operator spaces and Banach space Geometry \jour J.
Funct. Anal. \vol 109 \yr1992 \pages 113--129\endref
 
\ref \key Pa6  \bysame \paper   The maximal operator space of a normed space
\jour Proc. Edinburgh Math. Soc. \toappear\endref


 \ref \key PaS \by V. Paulsen and R. Smith \paper Multilinear maps and
tensor norms on operator systems \jour J. Funct. Anal. \vol 73
\yr1987 \pages 258--276\endref


 


\ref \key Pel \by V. Peller \paper  Estimates of functions of power
bounded operators on Hilbert space \jour J. Oper. Theory  \vol 7
\yr1982 \pages 341-372\endref

 \ref \key Pes \by V. Pestov \paper Operator spaces and residually
finite-dimensional $C^*$-algebras \jour
J. Funct. Anal. \vol 123 \yr1994 \pages 308--317\endref


\ref \key Pi1 \by J.P. Pier \book Amenable locally compact
groups \publ   Wiley Interscience \publaddr  New York \yr1984\endref

 \ref \key Pi2 \bysame \book Amenable Banach algebras  \publ Pitman, Longman
\yr1988\endref

\ref \key P1  \by G. Pisier \book Similarity problems and
completely bounded maps \publ Springer Lecture Notes \vol 1618 \yr1995\endref


\ref\key P2   \bysame \paper Multipliers and lacunary sets
in non amenable groups 
\jour Amer. J. Math. \vol 117  \yr1995 \pages 337--376\endref
 
 
\ref\key P3  \bysame \paper An introduction to the theory of
operator spaces \paperinfo
Preprint \toappear\endref

 \ref \key P4  \bysame \paper  A simple proof
of a theorem of Kirchberg
and related results on
$C^*$-norms \jour
J. Op. Theory. \vol 35 \yr1996 \pages 317--335\endref

\ref \key P5   \bysame \paper Un op\'erateur 
polyn{\^o}mialement born\'e sur un Hilbert qui n'est pas semblable
\`a une contraction \jour Comptes Rendus Acad. Sci. Paris \vol 322 
\yr1996 
\pages 547-550 \paperinfo (and article to appear
in  Journal Amer. Math. Soc.)\endref

\ref \key P6 \bysame \paper The operator Hilbert space $OH$,
complex interpolation and tensor norms \jour
Memoirs Amer.
Math. Soc.  \vol 122, 585 \yr1996 \pages 1--103\endref
 
 

\ref\key R \by  Z.J. Ruan \paper Subspaces of $C^*$-algebras \jour
J. Funct. Anal. \vol 76
 \yr1988 217--230\endref
 
\ref \key Ru\by  W. Rudin \book Fourier analysis on groups \publ
Interscience \publaddr New-York \yr1962\endref
  
\ref\key Sh \by M.V. Sheinberg \paper   A characterization of the algebra 
$C(\Omega)$
in terms of cohomology groups \jour Uspekhi Mat. Nauk \vol  32 \yr1977 \pages 
203--204 \transl\nofrills Russian\endref

 \ref \key SN \by  B.  Sz.-Nagy  \paper On uniformly bounded linear
transformations on Hilbert space \jour
Acta Sci. Math.  (Szeged)  \vol 11  \yr1946-48 \pages 152--157\endref


 \ref\key  V \by N. Varopoulos \paper A theorem on operator algebras \jour Math. 
Scand \vol
37 \yr1975 \pages 173--182\endref
 

\ref\key W \by S. Wassermann \paper
 On tensor products of certain
group $C^*$-algebras \jour J. Funct. Anal. \vol 23 \yr1976 \pages 
239--254\endref
 
\endRefs
\enddocument